\documentclass[preprint,3p]{elsarticle}
\usepackage{cite}
\usepackage{amsmath}
\usepackage{graphicx,wrapfig}
\usepackage{amssymb,amsgen,color}
\usepackage{amsfonts}
\usepackage{hyperref}
\usepackage{multirow}
\usepackage{epsfig,epstopdf,color,bm}

\usepackage{braket,amsfonts}
\usepackage{booktabs}
\newcommand{\be}{\begin{equation}}
\newcommand{\ee}{\end{equation}}

\newcommand{\eps}{\varepsilon}
\newtheorem{theorem}{Theorem}[section]
\newtheorem{lemma}{Lemma}[section]
\newtheorem{remark}{Remark}[section]
\newtheorem{corollary}{Corollary}[section]

\biboptions{numbers,sort&compress}
\makeatletter 
\@addtoreset{equation}{section}
\@addtoreset{figure}{section}
\@addtoreset{table}{section}
\makeatother  

\graphicspath{{figures/}}

\begin{document}
	\begin{frontmatter}
		
\title{Spatial resolution of different discretizations over long-time for the Dirac equation with small potentials}

\author[mymainaddress]{Yue Feng}
\ead{fengyue@u.nus.edu}

\author[mymainaddress,mysecondaryaddress]{Jia Yin\corref{mycorrespondingauthor}}
\cortext[mycorrespondingauthor]{Corresponding author}
\ead{jiayin@lbl.gov}

\address[mymainaddress]{Department of Mathematics, National University of Singapore, Singapore 119076, Singapore}
\address[mysecondaryaddress]{Computational Research Division, Lawrence Berkeley National Laboratory, Berkeley,\\ CA 94720, USA}
				
\begin{abstract}
We compare the long-time error bounds and spatial resolution of finite difference methods with different spatial discretizations for the Dirac equation with small electromagnetic potentials characterized by $\varepsilon \in (0, 1]$ a dimensionless parameter. We begin with the simple and widely used finite difference time domain (FDTD) methods, and establish rigorous error bounds of them, which are valid up to the time at $O(1/\varepsilon)$. In the error estimates, we pay particular attention to how the errors depend explicitly on the mesh size $h$ and time step $\tau$ as well as the small parameter $\varepsilon$. Based on the results, in order to obtain ``correct'' numerical solutions up to the time at $O(1/\varepsilon)$, the $\varepsilon$-scalability (or meshing strategy requirement) of the FDTD methods should be taken as $h = O(\varepsilon^{1/2})$ and $\tau = O(\varepsilon^{1/2})$. To improve the spatial resolution capacity, we apply the Fourier spectral method to discretize the Dirac equation in space. Error bounds of the resulting finite difference Fourier pseudospectral (FDFP) methods show that they exhibit uniform spatial errors in the long-time regime, which are optimal in space as suggested by the Shannon's sampling theorem.  Extensive numerical results are reported to confirm the error bounds and demonstrate that they are sharp.
\end{abstract}

		
\end{frontmatter}

\section{Introduction}
The Dirac equation is a relativistic wave equation derived by the British physicist Paul Dirac in 1928 \citep{Dirac1,Dirac2}. It plays a fundamental role in relativistic particle physics and is consistent with both the principles of quantum mechanics and the theory of special relativity. The Dirac equation describes the motion of the spin-1/2 massive particles such as electrons and positrons. It fully accounts for the details of the hydrogen spectrum and predicts the existence of a new form of matter, antimatter \citep{Anderson}. In the past decades, the Dirac equation has been widely used to study the structures and dynamical properties of graphene, graphite and other two-dimensional (2D) materials \citep{AMP,FW2,NGPNG,NGMJ}. Furthermore, some works are devoted to the analysis and numerical simulations of the Dirac points and the Shr\"odinger equation with the honeycomb lattice potential \citep{AZ,FW1}. There are also some applications in the relativistic effects of molecules in super intense lasers \citep{BSG,FLB}, and the motion of nucleons in the covariant density function theory in the relativistic framework \citep{BCJY, LMZ,Ring}.

We consider the following dimensionless Dirac equation on the unit torus in $d$ dimensions ($d = 1, 2, 3$) 
\begin{equation}
\label{eq:Dirac}
i\partial_t\Psi =  \left(- i\sum_{j = 1}^{d}
	\alpha_j\partial_j + \beta \right)\Psi+ \varepsilon\left(V(t, \mathbf{x})I_4 - \sum_{j = 1}^{d}A_j(t, \mathbf{x})\alpha_j\right)\Psi,
	\quad \mathbf{x}\in{\mathbb{T}^d}.
\end{equation}
Here, $i = \sqrt{-1}$, $0 < \varepsilon \leq 1$ is a dimensionless parameter, $t$ is time, $\mathbf{x} \in \mathbb{T}^d$ is the spatial coordinate vector, $\partial_j = \frac{\partial}{\partial x_j}( j = 1, 2, 3)$, $\Psi:=\Psi(t, \mathbf{x}) =  (\psi_1(t, \mathbf{x})$, $\psi_2(t, \mathbf{x}),
	\psi_3(t, \mathbf{x}),\psi_4(t, \mathbf{x}))^T \in \mathbb{C}^4$ is the complex-valued spinor wave function.  $I_n$ is the $n \times n$ identity matrix for $n \in \mathbb{N}$, $V:=V(t, \mathbf{x}) \in \mathbb{R}$ is the electric potential, and $\mathbf{A}:=\mathbf{A}(t, \mathbf{x}) = (A_1(t, \mathbf{x}), A_2(t, \mathbf{x}), A_3(t, \mathbf{x}))^T \in \mathbb{R}^3$ stands for the magnetic potential. In addition, the $4\times 4$ matrices $\alpha_1, \alpha_2, \alpha_3$ and $\beta$ are defined as
\begin{equation}
\label{matrices}
\alpha_1 = \begin{pmatrix} \mathbf{0} & \sigma_1\\ \sigma_1 & \mathbf{0}\end{pmatrix},\quad \alpha_2 = \begin{pmatrix} \mathbf{0} & \sigma_2\\ \sigma_2 & \mathbf{0}\end{pmatrix}, \quad \alpha_3 = \begin{pmatrix} \mathbf{0} & \sigma_3\\ \sigma_3 & \mathbf{0}\end{pmatrix}, \quad \beta = \begin{pmatrix} I_2 & \mathbf{0}\\ \mathbf{0} & -I_2\end{pmatrix},
\end{equation}
where $\sigma_j$ ($j=1,2,3$) are the Pauli matrices
	\begin{equation}
	\label{Pauli}
	\sigma_1 = \begin{pmatrix} 0 &\ \  1\\ 1 &\ \  0 \end{pmatrix}, \quad
	\sigma_2 = \begin{pmatrix} 0 & \ \ -i \\ i &\ \   0\end{pmatrix}, \quad
	\sigma_3 = \begin{pmatrix} 1 &\ \ 0 \\ 0 &\ \  -1 \end{pmatrix}.
	\end{equation}

Similar to that in \citep{BCJT}, in one dimension (1D) and two dimension (2D), the Dirac equation \eqref{eq:Dirac} can be reduced to the following two-component form with wave function $\Phi : = \Phi(t, {\mathbf{x}}) = (\phi_1(t, {\mathbf{x}}), \phi_2(t, {\mathbf{x}}))^T \in \mathbb{C}^2$ \citep{BCJY,BCY,BY}
\begin{equation}
\label{eq:Dirac_21}
i\partial_t\Phi =  \left(- i\sum_{j = 1}^{d}
	\sigma_j\partial_j + \sigma_3 \right)\Phi+ \varepsilon \left(V(t, \mathbf{x})I_2 - \sum_{j = 1}^{d}A_j(t, \mathbf{x})\sigma_j\right)\Phi,
\end{equation}
for $x \in \mathbb{T}^d (d=1,2)$, where $\Phi = (\psi_1, \psi_4)^T$ (or $\Phi = (\psi_2, \psi_3)^T$).

In order to study the dynamics of the Dirac equation \eqref{eq:Dirac_21}, the initial condition is taken as 
\begin{equation}
\label{eq:initial}
\Phi(t=0, \mathbf{x}) = \Phi_0(\mathbf{x}), \quad x \in \mathbb{T}^d.
\end{equation} 
The Dirac equation \eqref{eq:Dirac_21} with \eqref{eq:initial} is dispersive, time symmetric and conserves the total probability \citep{BCJT}
\begin{equation}
\|\Phi(t, \cdot)\|^2 : = \int_{\mathbb{T}^d} |\Phi(t, \mathbf{x})|^2 d \mathbf{x} = \int_{\mathbb{T}^d} \sum^2_{j=1} |\phi_j(t, \mathbf{x})|^2 d\mathbf{x} \equiv \|\Phi(0, \cdot)\|^2 =: \|\Phi_0\|^2, \quad t \geq 0.
\end{equation}
If the electromagnetic potentials are time-independent, i.e., $V(t, \mathbf{x}) = V(\mathbf{x})$ and $A_i(t, \mathbf{x}) = A_i(\mathbf{x})$ for $j=1, \cdots, d$, the following energy functional is also conserved
\begin{align}
E(\Phi(t, \cdot)) & : = \int_{\mathbb{T}^d} \left( -i \sum_{j=1}^d \Phi^{\ast}\sigma_j\partial_j \Phi + \Phi^{\ast}\sigma_3\Phi +\eps\left(V(\mathbf{x}) |\Phi|^2 - \sum_{j=1}^d A_j(\mathbf{x}) \Phi^{\ast}\sigma_j\Phi\right)\right) d\mathbf{x} \nonumber\\
	& \equiv E(\Phi(0, \cdot)) =: E(\Phi_0), \quad t \geq 0,
\end{align}
where $\bar{f}$ denotes the complex conjugate of $f$ and $\Phi^{\ast} = {\bar{\Phi}}^T$.

Introducing the total probability density $\rho : = \rho(t, \mathbf{x})$ as 
\begin{equation}
\label{eq:rho}
	\rho(t, \mathbf{x}) = \sum_{j=1}^2 \rho_j(t, \mathbf{x}) = \Phi(t, \mathbf{x})^{\ast}\Phi(t, \mathbf{x}), \quad \mathbf{x} \in \mathbb{T}^d,
\end{equation}
with the $j$th probability density $\rho_j := \rho_j(t, \mathbf{x})$ given as
\begin{equation}
 \rho_j(t, \mathbf{x})  = |\phi_j(t, \mathbf{x}) |^2,\quad 	\mathbf{x} \in \mathbb{T}^d,
\end{equation}
and the current density $\mathbf{J}(t, \mathbf{x}) = (J_1(t, \mathbf{x}), J_2(t, \mathbf{x}))^T$ with 
\begin{equation}
\label{eq:J}
J_l(t, \mathbf{x}) : = \Phi(t, \mathbf{x})^{\ast} \sigma_l \Phi(t, \mathbf{x}), \quad l = 1, 2, 	
\end{equation}
we can derive the following conservation law from the Dirac equation \citep{BCJT}
\begin{equation}
\partial_t \rho(t,\mathbf{x}) + \nabla \cdot \mathbf{J}(t, \mathbf{x}) = 0, \quad \mathbf{x} \in \mathbb{T}^d,\quad t \geq 0. 	
\end{equation}

 For the Dirac equation \eqref{eq:Dirac} (or \eqref{eq:Dirac_21}) with $\varepsilon = 1$, there are extensive analytical and numerical studies in the literature. Along the analytical front, for the existence and multiplicity of bound states and/or standing wave solutions, we refer to \citep{Das1,Das2,ES,GGT,Gross} and references therein. For the numerical aspects, different numerical methods have been presented and analyzed in the literature, such as the finite difference time domain (FDTD) methods \citep{BCJT, BHM,MY}, exponential wave integrator Fourier pseudospectral (EWI-FP) method \citep{BCJT}, time-splitting Fourier pseudospectral (TSFP) method \citep{BCY,BY,FLB}, WKB method \citep{CMS,SM1,Spohn} and Gaussian beam method \citep{JWY,WHJY}. For more details related to the numerical schemes, we refer to \citep{AH,BCJT2016,BSG,FLB,FLB1,FLB2,Gosse,GSX,Shebalin,XST} and references therein. 

However, to the best of our knowledge, so far there are few numerical analysis results on the error bounds of different numerical methods for the long-time dynamics of the Dirac equation \eqref{eq:Dirac} in the literature, especially the error bounds which are valid up to the time at $T_{\varepsilon} = O(1/\varepsilon)$, and their explicit dependence on the mesh size $h$, time step $\tau$ as well as the small parameter $\varepsilon \in (0, 1]$. Formally, by the energy estimate and Gronwall inequality, error bounds in the finite time for $t \in [0,T]$ behave like $O(e^{CT}\tau^{p}) $ with $p$th order temporal discretization, which indicates the unbounded temporal error bound $O(e^{CT/\varepsilon}\tau^{p}) $ for $t \in [0, T/\varepsilon]$. In order to carry out valid error estimates in the long-time regime, we begin with the proper setup, i.e., the Dirac equation with $O(\varepsilon)$ potential. In our previous work, the long-time error bounds of different numerical schemes for the nonlinear Klein-Gordon equation with weak nonlinearity have been rigorously established \citep{BFY,Feng}. A similar problem for the Schr\"{o}dinger equation by using splitting methods has been studied \citep{Dujardin1,Dujardin2}. The aim of this paper is to carry out rigorous error bounds of finite difference methods with different spatial discretizations for the Dirac equation \eqref{eq:Dirac} in the long-time regime. We begin with the detailed analysis on the stability and convergence of four frequently used finite difference time domain (FDTD) methods including the implicit/semi-implicit/explicit scheme. Error bounds of the FDTD methods indicate that they are \textbf{under-resolution} both in space and time with respect to $\varepsilon \in (0, 1]$ regarding to the Shannon's sampling theorem - to resolve a wave one needs a few points per wavelength, since the wavelength in space and time is at $O(1)$. Then, in order to improve the spatial resolution capacity of the FDTD methods, we adapt the Fourier spectral method to discretize the Dirac equation \eqref{eq:Dirac} in space, and rigorous error estimates of finite difference Fourier pseudospectral (FDFP) methods show that the spatial errors of FDFP methods are uniform for the long-time dynamics of the Dirac equation. Thus, FDFP methods offer compelling advantage over FDTD methods in spatial resolution when $0 < \varepsilon \ll 1$.

The rest of the paper is organized as follows. In Section 2, we adapt the finite difference methods to discretize the Dirac equation in time, combined with the centered finite difference method and Fourier spectral method in spatial discretization. In Section 3, we give the conservation properties, and analyze the stability conditions of the FDTD and FDFP methods. In Section 4, the long-time error bounds of the FDTD and FDFP methods for the Dirac equation up to the time at $O(1/\varepsilon)$ are rigorously presented. Extensive numerical results are reported in Section 5 to confirm our error bounds. Finally, some conclusions are drawn in Section 6. Throughout this paper, we adopt the notation $A \lesssim B$ to represent that there exists a generic constant $C > 0$, which is independent of the mesh size $h$ and time step $\tau$ as well as $\varepsilon$ such that $|A| \leq C B$.

\section{Finite difference methods with different spatial discretizations}
In this section, we use the finite difference methods,  including the Crank-Nicolson finite difference method, two semi-implicit finite difference methods and the leap-frog finite difference method, to discretize the Dirac equation \eqref{eq:Dirac_21} in time. For spatial discretization, we adapt the centered finite difference method and the Fourier spectral method. For simplicity of notations, we only present the numerical methods and their analysis for the Dirac equation \eqref{eq:Dirac_21} in 1D. Generalization to \eqref{eq:Dirac} and/or higher dimensions is straightforward and results are valid without modifications. In 1D, the Dirac equation \eqref{eq:Dirac_21} with periodic boundary conditions collapses to 
\begin{align}
\label{eq:Dirac_1D}
&i\partial_t\Phi =  (- i \sigma_1 \partial_x + \sigma_3 )\Phi+ \varepsilon(V(t, x)I_2 - A_1(t, x)\sigma_1)\Phi, \quad x \in \Omega = (a, b), \quad t > 0,\\
\label{eq:ib}
&\Phi(t, a) = \Phi(t, b),\quad \partial_x \Phi(t, a) = \partial_x \Phi(t, b),\quad t \geq 0; \quad \Phi(0, x) = \Phi_0(x),\quad x \in \overline{\Omega}, 
\end{align}
where $\Phi := \Phi(t, x)$, $\Phi_0(a) = \Phi_0(b)$ and $\Phi_0'(a) = \Phi_0'(b)$.

\subsection{Semi-discretization by finite difference in time}
Take time step size $\tau = \Delta t > 0$ and denote the time steps as 
\begin{equation}
t_n :=n\tau, \quad n = 0,1,2,\cdots.
\end{equation}
For $n\geq 0$, denote $\Phi^n(x)$ to be the semi-discretized numerical approximation of $\Phi(t_n, x)$, $V^n(x) = V(t_n, x)$, $V^{n+1/2}(x) = V(t_n+\tau/2, x)$, $A_{1}^n(x) = A_1(t_n, x)$, and $A_{1}^{n+1/2}(x) = A_1(t_n+\tau/2, x)$. Introduce the following finite difference discretization operators in time:
\begin{equation*}
\delta^+_t \Phi^n(x) = \frac{\Phi^{n+1}(x) - \Phi^n(x)}{\tau},\quad \delta_t \Phi^n(x) = \frac{\Phi^{n+1}(x) - \Phi^{n-1}(x)}{2\tau},
\end{equation*}
and introduce the average vector:
\begin{equation*}
\Phi^{n+\frac{1}{2}}(x)	 = \frac{\Phi^{n+1}(x) + \Phi^n(x)}{2},
\end{equation*}
then four frequently used finite difference discretization in time for the Dirac equation \eqref{eq:Dirac_1D} could be expressed as:

I. Crank-Nicolson method
\begin{equation}
\label{eq:CN}
i\delta^+_t \Phi^n(x) = \left(-i\sigma_1\frac{d}{dx}+ \sigma_3\right)\Phi^{n+1/2}(x) + \varepsilon \left(V^{n+1/2}(x) I_2 -A^{n+1/2}_{1}(x)\sigma_1\right)\Phi^{n+1/2}(x),\quad  n \geq 0.	
\end{equation}

II. A semi-implicit method
\begin{equation}
\label{eq:SI1}
i\delta_t \Phi^n(x) = -i\sigma_1\frac{d}{dx} \Phi^n(x) +\left(\sigma_3 +  \varepsilon (V^n(x) I_2 -A^n_{1}(x)\sigma_1)\right)\frac{\Phi^{n+1}(x)+\Phi^{n-1}(x)}{2},\quad  n \geq 1.	
\end{equation}

III. Another semi-implicit method
\begin{equation}
\label{eq:SI2}
i\delta_t \Phi^n(x) = \left(-i\sigma_1\frac{d}{dx}+ \sigma_3\right)\frac{\Phi^{n+1}(x)+ \Phi^{n-1}(x)}{2} + \varepsilon\left(V^n(x) I_2 -A^n_{1}(x)\sigma_1\right)\Phi^{n}(x),\quad n \geq 1.	
\end{equation}

IV. Leap-frog method
\begin{equation}
\label{eq:LF}
i\delta_t \Phi^n(x) = \left(-i\sigma_1\frac{d}{dx} + \sigma_3\right)\Phi^n(x) +\varepsilon \left(V^n(x) I_2 -A^n_{1}(x)\sigma_1\right)\Phi^n(x),\quad n \geq 1.	
\end{equation}

The initial and boundary conditions for these methods are taken as
\begin{equation}
\label{eq:init}
\Phi^n(a) = \Phi^n(b), \quad \frac{d}{dx}\Phi^n(a) = \frac{d}{dx}\Phi^n(b), \quad n\geq 0; \quad \Phi^0(x) = \Phi_0(x), \quad x\in (a, b), 
\end{equation}
and the first step for \eqref{eq:SI1}, \eqref{eq:SI2} and \eqref{eq:LF} can be updated by using Taylor expansion and noticing \eqref{eq:Dirac_1D} as 
\begin{equation}
\label{eq:1st}
\Phi^1(x) = \Phi_0(x) -\tau\sigma_1\Phi'_0(x) - i\tau\left(\sigma_3 + \varepsilon V^0(x) I_2-\varepsilon A^0_{1}(x)\sigma_1\right)\Phi_0(x), \quad  x\in (a, b).	
\end{equation}

\subsection{Spatial discretization by finite difference method}
We first use finite difference method to discretize \eqref{eq:CN} to \eqref{eq:LF} in space. Choose mesh size $h : = \Delta x = (b - a)/N$ with $N$ being an even positive integer, and denote the grid points
\begin{equation}
x_j :=a +j h,\quad j = 0, 1, \cdots, N.
\end{equation}
Take $X_N = \{U=(U_0, U_1,\cdots,U_N)^T | U_j \in \mathbb{C}^2, j = 0,1, \cdots, N, U_0 = U_N\}$ and let $U_{-1} = U_{N-1}$, $U_{N+1} = U_1$ if they are involved. Define the index set $
\mathcal{T}_N = \{l | l = -N/2,-N/2+1,\cdots,N/2-1\}$, and $\mu_l = \frac{2\pi l}{b-a}$ for $l \in \mathcal{T}_N$,
then for any $U \in X_N$, its Fourier representation can be expressed as 
\begin{equation}
U_j = \sum_{l \in \mathcal{T}_N}\widetilde{U}_l e^{i\mu_l(x_j-a)} = \sum_{l \in \mathcal{T}_N}\widetilde{U}_l e^{2ijl\pi/N}, \quad j = 0,1, \cdots, N,
\end{equation}
where 
\begin{equation}
\widetilde{U}_l = \frac{1}{N} \sum_{j=0}^{N-1} U_j e^{-2ijl\pi/N}, \quad l \in \mathcal{T}_N.
\end{equation}
The standard $l^1$ and $l^2$-norms in the space $X_N$ are given as
\begin{equation}
\|U\|_{l^1} = h \sum_{j=0}^{N-1} |U_j|,\quad \|U\|^2_{l^2} = h \sum_{j=0}^{N-1} |U_j|^2,\quad U \in X_N.
\end{equation}
For $0 \leq j \leq N$ and $n \geq 0$, let $\Phi^n_j$ be the full-discretized numerical approximation of $\Phi(t_n, x_j)$, $V^n_j = V(t_n, x_j)$, $V^{n+1/2}_j = V(t_n+\tau/2, x_j)$, $A^n_{1, j} = A_1(t_n, x_j)$ and $A^{n+1/2}_{1, j} = A_1(t_n + \tau/2, x_j)$. Denote $\Phi^n = (\Phi^n_0, \Phi^n_1,\cdots, \Phi^n_N)^T \in X_N$ as the solution at $t=t_n$. Introduce the following finite difference discretization operators
\begin{equation*}
\delta^+_t \Phi^n_j = \frac{\Phi^{n+1}_j - \Phi^n_j}{\tau},\quad \delta_t \Phi^n_j = \frac{\Phi^{n+1}_j - \Phi^{n-1}_j}{2\tau},\quad \delta_x \Phi^n_j = \frac{\Phi^n_{j+1} - \Phi^n_{j-1}}{2h},
\end{equation*}
and the average vector
\begin{equation*}
\Phi^{n+\frac{1}{2}}_j	 = \frac{\Phi^{n+1}_j + \Phi^n_j}{2},
\end{equation*}
then by discretizing \eqref{eq:CN} to \eqref{eq:LF} in space, we have the following four finite difference time domain (FDTD) methods for $j = 0,1,\cdots, N-1$:

I. Crank-Nicolson finite difference (CNFD) method
\begin{equation}
\label{eq:CNFD}
i\delta^+_t \Phi^n_j = \left(-i\sigma_1\delta_x+ \sigma_3\right)\Phi^{n+1/2}_j + \varepsilon \left(V^{n+1/2}_j I_2 -A^{n+1/2}_{1, j}\sigma_1\right)\Phi^{n+1/2}_j,\quad n \geq 0.	
\end{equation}

II. A semi-implicit finite difference (SIFD1) method
\begin{equation}
\label{eq:SIFD1}
i\delta_t \Phi^n_j = -i\sigma_1\delta_x \Phi^n_j +\left(\sigma_3 +  \varepsilon (V^n_j I_2 -A^n_{1, j}\sigma_1)\right)\frac{\Phi^{n+1}_j+\Phi^{n-1}_j}{2},\quad n \geq 1.	
\end{equation}

III. Another semi-implicit finite difference (SIFD2) method
\begin{equation}
\label{eq:SIFD2}
i\delta_t \Phi^n_j = \left(-i\sigma_1\delta_x+ \sigma_3\right)\frac{\Phi^{n+1}_j+ \Phi^{n-1}_j  }{2} + \varepsilon\left(V^n_j I_2 -A^n_{1, j}\sigma_1\right)\Phi^{n}_j,\quad n \geq 1.	
\end{equation}

IV. Leap-frog finite difference (LFFD) method
\begin{equation}
\label{eq:LFFD}
i\delta_t \Phi^n_j = \left(-i\sigma_1\delta_x + \sigma_3\right)\Phi^n_j +\varepsilon \left(V^n_j I_2 -A^n_{1, j}\sigma_1\right)\Phi^n_j,\quad n \geq 1.	
\end{equation}

The initial and boundary conditions in \eqref{eq:ib} for the FDTD methods are discretized as
\begin{equation}
\label{eq:she0}
\Phi^{n+1}_0 = \Phi^{n+1}_N,\quad	 \Phi^{n+1}_{-1} = \Phi^{n+1}_{N-1}, \quad n \geq 0;\quad \Phi^0_j = \Phi_0(x_j), \quad j = 0,1,\cdots, N,
\end{equation}
and the first step for SIFD1 \eqref{eq:SIFD1}, SIFD2 \eqref{eq:SIFD2} and LFFD \eqref{eq:LFFD} could be directly derived from \eqref{eq:1st} as
\begin{equation}
\label{eq:phi1}
\Phi^1_j = \Phi_0(x_j) -\tau\sigma_1\Phi'_0(x_j) - i\tau\left(\sigma_3 + \varepsilon V^0_j I_2-\varepsilon A^0_{1, j}\sigma_1\right)\Phi_0(x_j), \quad j =0,1, \cdots, N.	
\end{equation}

The above four FDTD methods are all time symmetric, i.e., they are unchanged under $\tau \leftrightarrow -\tau$ and $n+1 \leftrightarrow n-1$ in the SIFD1, SIFD2 and LFFD methods or $n+1 \leftrightarrow n$ in the CNFD method. Besides, their memory costs are all $O(N)$. The CNFD method \eqref{eq:CNFD} is implicit and at each time step for $n \geq 0$, the corresponding linear system is coupled so that it needs to be solved through either a direct solver or an iterative solver with the computational cost per time step depending on the solver, which is usually much larger than $O(N)$, especially in 2D and 3D. The SIFD1 method \eqref{eq:SIFD1} is an implicit scheme, however, the corresponding linear system is decoupled and can be solved explicitly at each time step as
\begin{equation}
\Phi^{n+1}_j = \left[(i-\varepsilon\tau V^n_j)I_2 - \tau\sigma_3 + \varepsilon\tau A^n_{1, j} \sigma_1\right]^{-1}H^n_j, \quad j =0 ,1, \cdots, N-1,
\end{equation}
with $H^n_j = -2i\tau \sigma_1\delta_x \Phi^n_j +\left[\left(i+\varepsilon\tau V^n_j\right)I_2+\tau\left(\sigma_3- \varepsilon A^n_{1, j} \sigma_1\right)\right]\Phi^{n-1}_j$, and its computational cost per time step is also $O(N)$. The SIFD2 method \eqref{eq:SIFD2} is also implicit, but the corresponding linear system can be decoupled in the phase (Fourier) space  and can be solved explicitly as
\begin{equation}
\widetilde{(\Phi^{n+1})}_l = \left(i I_2 - \frac{\tau\sin(\mu_l h)}{h}\sigma_1-\tau\sigma_3\right)^{-1} L^n_l,\quad l \in \mathcal{T}_N,	
\end{equation}
where
\begin{equation*}
L^n_l = \left(i I_2 +  \frac{\tau\sin(\mu_l h)}{h}\sigma_1 + \tau\sigma_3\right) \widetilde{(\Phi^{n-1})}_l + 2\tau \widetilde{(G^n\Phi^n)}_l,
\end{equation*}
and $G^n = (G^n_0, G^n_1, \cdots, G^n_N)^T \in X_N$ with $G^n_j = \varepsilon(V^n_j I_2 -A^n_{1, j}\sigma_1)$ for $j = 0,1,\cdots, N$. Therefore, its computational cost per time step is $O(N \ln N)$. The LFFD method \eqref{eq:LFFD} is an explicit scheme and  might be the simplest and most efficient numerical method for the Dirac equation with the computational cost per time step at $O(N)$.  Based on the analysis on the computational cost, the LFFD method is the most efficient numerical method and the CNFD methods is the most expensive one.

\subsection{Spatial discretization by Fourier spectral method}
We could also use Fourier spectral method to discretize \eqref{eq:CN} to \eqref{eq:LF} in space. Denote
\begin{equation*}
	Y_N = Z_N \times Z_N, \quad Z_N = \mbox{span}\{\phi_l(x) = e^{i\mu_l(x-a)},\ l \in \mathcal{T}_N\}
\end{equation*}
Let $(C_p(\overline{\Omega}))^2$ be the function space consisting of all periodic vector function $U(x): \overline{\Omega} = [a, b] \to \mathbb{C}^2$. For any $U(x) \in (C_p(\overline{\Omega}))^2$ and $U \in X_N$, define $P_N: (L^2(\Omega))^2 \to Y_N$ as the standard projector \citep{ST}, $I_N: (C_p(\overline{\Omega}))^2 \to Y_N$ and $I_N : X_{N} \to Y_N$ as the standard interpolation operator, i.e., for any $a \leq x \leq b$
\begin{equation}
	(P_N U)(x) = \sum_{l \in \mathcal{T}_N}\widehat{U}_l e^{i\mu_l(x-a)}, \quad (I_N U)(x) = \sum_{l \in \mathcal{T}_N}\widetilde{U}_l e^{i\mu_l(x-a)},
\end{equation} 
with 
\begin{equation}
\label{eq:Fc}
\widehat{U}_l = \frac{1}{b-a}\int^b_a U(x) e^{-i\mu_l(x-a)} dx,\quad \widetilde{U}_l = \frac{1}{N}\sum_{j=0}^{N-1}U_j e^{-2ijl\pi/N},\quad l \in \mathcal{T}_N, 	
\end{equation}
where $U_j = U(x_j)$ when $U$ is a function.

The Fourier spectral discretization for \eqref{eq:CN} to \eqref{eq:LF} is to project the equations onto $Y_N$ by using the projection operator $P_N$. Let $\Phi^n(x)$ be the approximation of $\Phi(t_n, x) (n \geq 0)$. Introduce the spectral differential operator $D^f_x$ to approximate the operator $\frac{d}{d x}$ as
\begin{equation}
D^f_x \Phi |_{x = x_j} = i \sum_{l \in \mathcal{T}_N} \mu_l \widehat{\Phi}_l e^{i\mu_l(x_j - a)},	
\end{equation}
and the average function
\begin{equation*}
\Phi^{n+\frac{1}{2}}	 = \frac{\Phi^{n+1} + \Phi^n}{2},	
\end{equation*}
then we have the following four finite difference Fourier spectral (FDFS) methods:

I. Crank-Nicolson Fourier spectral (CNFS) method
\begin{equation}
\label{eq:CNFS}
i\delta^+_t \Phi^n_j = -i\sigma_1 D^f_x \Phi^{n+1/2}|_{x = x_j}+ \sigma_3\Phi^{n+1/2}_j + \varepsilon \left(V^{n+1/2}_j I_2 -A^{n+1/2}_{1, j}\sigma_1\right)\Phi^{n+1/2}_j,\quad n \geq 0.	
\end{equation}

II. A semi-implicit Fourier spectral (SIFS1) method
\begin{equation}
\label{eq:SIFS1}
i\delta_t \Phi^n_j = -i\sigma_1 D^f_x  \Phi^n|_{x = x_j} +\left(\sigma_3 +  \varepsilon (V^n_j I_2 -A^n_{1, j}\sigma_1)\right)\frac{\Phi^{n+1}_j+\Phi^{n-1}_j}{2},\quad n \geq 1.	
\end{equation}

III. Another semi-implicit Fourier spectral (SIFS2) method
\begin{equation}
\label{eq:SIFS2}
i\delta_t \Phi^n_j = -i\sigma_1D^f_x \frac{\Phi^{n+1}+ \Phi^{n-1}}{2}|_{x = x_j} + \sigma_3\frac{\Phi^{n+1}_j+ \Phi^{n-1}_j}{2} + \varepsilon\left(V^n_j I_2 -A^n_{1, j}\sigma_1\right)\Phi^{n}_j,\quad n \geq 1.	
\end{equation}

IV. Leap-frog Fourier spectral (LFFS) method
\begin{equation}
\label{eq:LFFS}
i\delta_t \Phi^n_j = -i\sigma_1 D^f_x \Phi^n|_{x = x_j} + \sigma_3\Phi^n_j +\varepsilon \left(V^n_j I_2 -A^n_{1, j}\sigma_1\right)\Phi^n_j,\quad n \geq 1.
\end{equation}
The initial condition is taken to be
\begin{equation}
\Phi_j^0 = (P_N\Phi_0)(x_j), \quad j = 0, 1, ..., N,
\label{eq:LFFS_initial}
\end{equation}
and the first step for SIFS1 \eqref{eq:SIFS1}, SIFS2 \eqref{eq:SIFS2}, and LFFS \eqref{eq:LFFS} is the same as that for FDTD methods.

In fact, in practice, it is difficult to compute the Fourier coefficients through integrals in \eqref{eq:Fc}. Therefore, we replace the projections by interpolations in numerical computation, which refers to the Fourier pseudospectral discretization. The resulting methods are called finite difference Fourier pseudospectral (FDFP) methods, including CNFP, SIFP1, SIFP2, and LFFP. Choosing $\Phi^0_j = \Phi_0(x_j) (j = 0, 1, \cdots, N)$, the FDFP methods for computing $\Phi_j^{n}$ ($n\geq 0$) are given as
\begin{equation}
\label{eq:FP}
\Phi^{n}_j = \sum_{l \in \mathcal{T}_N} \widetilde{(\Phi^{n})}_l e^{2ijl\pi/N},\quad j = 0, 1, \cdots, N.
\end{equation}
The formulas of the FDFP methods could be expressed in the same form as \eqref{eq:CNFS} to \eqref{eq:LFFS}, with
\begin{align}
  D_x^f\Phi|_{x = x_j} = & \ i\sum_{l\in\mathcal{T}_N}\frac{\mu_l}{b-a}\int_a^b\Phi(x)e^{-i\mu_l(x-a)}dx\ e^{i\mu_l(x_j-a)}\nonumber\\
  = & \ \frac{i}{b-a}\int_a^b\sum_{l\in\mathcal{T}_N}\mu_le^{i\mu_l(x_j-x)}\Phi(x)dx \label{eq:D_FS}
\end{align}
for FDFS methods, and
\begin{align}\label{eq:D_FP}
D_x^f\Phi|_{x = x_j} = & \ \frac{i}{N}\sum_{k=0}^{N-1}\sum_{l\in\mathcal{T}_N}\mu_le^{2i(j-k)l\pi/N}\Phi^n_k
\end{align}
for FDFP methods.

\section{Conservation and linear stability}
The Crank-Nicolson methods conserve mass and energy. Here we give proof for the CNFP method, and proof for the CNFD method is similar. Besides, stability conditions of FDTD and FDFP methods are derived.

\subsection{Mass and energy conservation for CNFD and CNFP}
For the CNFD method \eqref{eq:CNFD}, we have the following conservative properties.
\begin{lemma}\label{lemma:conserv1}
The CNFD method \eqref{eq:CNFD} conserves the mass in the discretized level, i.e.,
\begin{equation}
\label{eq:mass}
\|\Phi^n\|^2_{l^2} := h \sum^{N-1}_{j=0} |\Phi^n_j|^2 \equiv h \sum^{N-1}_{j=0} |\Phi^0_j|^2 = \|\Phi^0\|^2_{l^2}.
\end{equation}
Furthermore, if $V(t, x) = V(x)$ and $A_1(t, x) = A_1(x)$ are time-independent, then the CNFD method \eqref{eq:CNFD} conserves the discretized energy as well, that is
\begin{align}
\label{eq:energy_con}
E^n_h := \ & h\sum^{N-1}_{j=0} \left[-i (\Phi^n_j)^{\ast}\sigma_1\delta_x \Phi^n_j +(\Phi^n_j)^{\ast}\sigma_3 \Phi^n_j + \varepsilon V_j|\Phi^n_j|^2 - \varepsilon A_{1, j}(\Phi^n_j)^{\ast}\sigma_1 \Phi^n_j \right] \nonumber\\
\equiv \ & E^0_h, \quad n \geq 0,
\end{align}
where $V_j = V(x_j)$ and $A_{1, j} = A_1(x_j)$ for $j = 0,1, \cdots, N$.
\end{lemma}

For the CNFP method, we also have the following conservative properties.
\begin{lemma}
\label{lemma:conserv2}
The CNFP method  conserves the mass as
\begin{align}
\label{eq:mass2}
\|\Phi^n\|_{l^2} = \|\Phi^0\|_{l^2}.
\end{align}
Furthermore, if $V(t, x) = V(x)$ and $A_1(t, x) = A_1(x)$ are time-independent, then the CNFP method conserves the energy as well, that is
\begin{align}
\label{eq:energy_con2}
\mathcal{E}^n_h := \ & h\sum^{N-1}_{j=0} \left[-i (\Phi^n_j)^{\ast}\sigma_1D_x^f \Phi^n|_{x=x_j} +(\Phi^n_j)^{\ast}\sigma_3 \Phi^n_j + \varepsilon V_j|\Phi^n_j|^2 - \varepsilon A_{1, j}(\Phi^n_j)^{\ast}\sigma_1 \Phi^n_j \right] \nonumber\\
\equiv \ & \mathcal{E}^0_h, \quad n \geq 0.
\end{align}
\end{lemma}

\emph{Proof.}
	For the mass conservation, similar to the proof of the CNFD method, multiplying both sides of \eqref{eq:CNFS} from left by $h\tau (\Phi^{n+1/2}_j)^{\ast}$ and taking the imaginary part, we have for $j = 0, 1, \cdots, N-1$
\begin{equation}\label{eq:mc_FS}
h|\Phi^{n+1}_j|^2 = h|\Phi^n_j|^2 - \tau h\left[(\Phi^{n+1/2}_j)^{\ast}\sigma_1 D^f_x \Phi^{n+1/2}|_{x=x_j} + (D_x^f\Phi^{n+1/2})^*|_{x=x_j}\sigma_1\Phi_j^{n+1/2}\right].	
\end{equation}
Summing the above equation for $j = 0, 1, \cdots, N-1$ and combing the expression of the operator $D^f_x$ in \eqref{eq:D_FP}, we could get
 \begin{align}
 \|\Phi^{n+1}\|^2_{l^2} = &\ \|\Phi^{n}\|^2_{l^2} - \tau h\sum^{N-1}_{j=0}\Bigg[(\Phi^{n+1/2}_j)^{\ast}\sigma_1 \left(\frac{i}{N}\sum_{k=0}^{N-1}\sum_{l\in\mathcal{T}_N}\mu_le^{2i(j-k)l\pi/N}\Phi_k^{n+1/2}\right)\nonumber\\
 	&\; + \left(-\frac{i}{N}\sum_{k=0}^{N-1}\sum_{l\in\mathcal{T}_N}\mu_le^{-2i(j-k)l\pi/N}(\Phi_k^{n+1/2})^*\right)\sigma_1\Phi_j^{n+1/2}\Bigg] \nonumber\\
 = & \  \|\Phi^n\|^2_{l^2} - \frac{i\tau h}{N}\sum_{j=0}^{N-1}\sum_{k=0}^{N-1}\sum_{l\in\mathcal{T}_N}\mu_l\left[e^\frac{2i(j-k)l\pi}{N}(\Phi^{n+1/2}_j)^*\sigma_1\Phi^{n+1/2}_k - e^\frac{2i(k-j)l\pi}{N}(\Phi^{n+1/2}_k)^*\sigma_1\Phi^{n+1/2}_j\right] \nonumber \\
 = & \ \|\Phi^{n}\|^2_{l^2}, \quad n \geq 0,
 \label{al:massproof}
\end{align}
 which implies {\eqref{eq:mass2}} by induction.
	
For the energy conservation, multiplying both sides of \eqref{eq:CNFS} from left by $2h (\Phi^{n+1}_j - \Phi^n_j)^{\ast}$ and taking the real part, we obtain 
 \begin{align}
 & -h \text{Re}\left[i(\Phi^{n+1}_j - \Phi^n_j)^{\ast} \sigma_1 D^f_x(\Phi^{n+1}+\Phi^n)|_{x=x_j}\right] + h\left[(\Phi^{n+1}_j)^{\ast}\sigma_3 \Phi^{n+1}_j - (\Phi^{n}_j)^{\ast}\sigma_3 \Phi^{n}_j\right] \nonumber \\  
 & \quad+ \varepsilon h V_j (|\Phi^{n+1}_j|^2-|\Phi^{n}_j|^2) - \varepsilon h A_{1, j} \left[(\Phi^{n+1}_j)^{\ast}\sigma_1 \Phi^{n+1}_j - (\Phi^{n}_j)^{\ast}\sigma_1 \Phi^{n}_j\right] = 0.
 \end{align}
Summing the above equations for $ j = 0, 1, \cdots, N-1$ and noticing the summation by parts formula, the mass conservation, and the expression of $D_x^f$ \eqref{eq:D_FP}, we could get
\begin{equation*}
\begin{split}
&h \sum^{N-1}_{j=0} \text{Re} \left(i(\Phi^{n+1}_j - \Phi^n_j)^{\ast}\sigma_1 D^f_x (\Phi^{n+1} +\Phi^n)|_{x=x_j}\right) \\
=&\  ih\sum^{N-1}_{j=0}(\Phi^{n+1}_j)^{\ast}\sigma_1 D^f_x\Phi^{n+1}|_{x=x_j} - ih\sum^{N-1}_{j=0}(\Phi^{n}_j)^{\ast}\sigma_1 D^f_x\Phi^{n}|_{x=x_j},	
\end{split}
\end{equation*}
and
\begin{equation*}
\begin{split}
& - ih \sum^{N-1}_{j=0}(\Phi^{n+1}_j)^{\ast}\sigma_1 D^f_x \Phi^{n+1}|_{x=x_j} + ih\sum^{N-1}_{j=0}(\Phi^{n}_j)^{\ast}\sigma_1 D^f_x\Phi^{n}|_{x=x_j}  + h\sum^{N-1}_{j=0}\left((\Phi^{n+1}_j)^{\ast}\sigma_3 \Phi^{n+1}_j - (\Phi^n_j)^{\ast}\sigma_3\Phi^n_j\right) \\
& + \varepsilon h \sum^{N-1}_{j=0} V_j (|\Phi^{n+1}_j|^2 - |\Phi^n_j|^2)  - \varepsilon h \sum^{N-1}_{j=0} A_{1, j}\left((\Phi^{n+1}_j)^{\ast} \sigma_1\Phi^{n+1}_j - (\Phi^n_j)^{\ast}\sigma_1\Phi^n_j\right) = 0,
\end{split}
\end{equation*}
which immediately implies \eqref{eq:energy_con2}.
\hfill $\square$ \bigskip

\subsection{Linear stability analysis}
Following the von Neumann linear stability analysis \citep{Smith}, we analyze the stability of FDTD and FDFP methods. The proof is similar to that in \citep{BCJY,MY} and we omit the details here for brevity.
\begin{lemma}
\label{lemma:sta1}
Denote $\Omega_T = [0, T] \times \Omega$. For the case where the electromagnetic potentials are not constants, we take
\begin{equation}\label{eq:max}
V_{\max} := \max_{(t, x) \in \overline{\Omega}_T}|V(t, x)|	, \quad A_{1, \max} := \max_{(t, x) \in \overline{\Omega}_T} |A_1(t, x)|,
\end{equation}
(i) The CNFD method \eqref{eq:CNFD} with \eqref{eq:she0} is unconditionally stable, i.e., it is stable for any $h > 0$, $\tau > 0$ and $0 < \varepsilon \leq 1$.\\
(ii) For any $0 < \varepsilon \leq 1$, the SIFD1 method \eqref{eq:SIFD1} with \eqref{eq:she0} is stable under the stability condition
\begin{equation}
\label{eq:st_SIFD1}
0 < \tau \leq h,\quad h > 0.	
\end{equation}
(iii) For any $0 < \varepsilon \leq 1$, the SIFD2 method \eqref{eq:SIFD2} with \eqref{eq:she0} is stable under the stability condition
\begin{equation}
\label{eq:st_SIFD2}
0 < \tau \leq \frac{1}{V_{\max}+ A_{1, \max}},\quad h > 0.	
\end{equation}
(iv) For any $0 < \varepsilon \leq 1$, the LFFD method \eqref{eq:LFFD} with \eqref{eq:she0} is stable under the stability condition
\begin{equation}
\label{eq:st_LFFD}
0 < \tau \leq \frac{h}{V_{\max}h+\sqrt{h^2+(1+h A_{1, \max})^2}},\quad h > 0.	
\end{equation}
\end{lemma}

Similarly to FDTD methods, the stability conditions for FDFS methods could be derived as follows:
\begin{lemma}\label{lemma:sta2}
	(i) The CNFS method \eqref{eq:CNFS} with \eqref{eq:she0} is unconditionally stable, i.e., it is stable for any $h > 0$, $\tau > 0$ and $0 < \varepsilon \leq 1$.	\\
	(ii) For any $0 < \varepsilon \leq 1$, the SIFS1 method \eqref{eq:SIFS1} with \eqref{eq:she0} is stable under the stability condition
	\begin{equation}
	\label{eq:st_SIFS1}
	0 < \tau \leq \frac{h}{\pi},\quad h > 0.	
	\end{equation}
	(iii) For any $0 < \varepsilon \leq 1$, the SIFS2 method \eqref{eq:SIFS2} with \eqref{eq:she0} is stable under the stability condition
	\begin{equation}
	\label{eq:st_SIFS2}
	0 < \tau \leq \frac{1}{V_{\max} + A_{1, \max}}, \quad h > 0.	
	\end{equation}
	(iv) For any $0 < \varepsilon \leq 1$, the LFFS method \eqref{eq:LFFS} with \eqref{eq:she0} is stable under the stability condition
	\begin{equation}
	\label{eq:st_LFFS}
	0 < \tau \leq \frac{h}{V_{\max}h+\sqrt{h^2+(\pi + h A_{1, \max})^2}},\quad h > 0,
	\end{equation}
with $V_{\max}$ and $A_{1, \max}$ defined in \eqref{eq:max}.
\end{lemma}


\section{Error estimates}
In this section, we propose and prove the long-time error estimates for FDTD and FDFP methods.

\subsection{Main results}
Motivated by the analysis of the Dirac equation \citep{AS,Das1,ES,FW1}, we assume that the exact solution of the Dirac equation \eqref{eq:Dirac_1D} up to the time $t= T_0/\varepsilon$ satisfies:
\begin{equation*}
\begin{split}
\Phi \in C^3([0, T_0/\varepsilon];(L^{\infty}(\Omega))^2) & \cap C^2([0, T_0/\varepsilon];(W^{1, \infty}_p(\Omega))^2)  \cap C^1([0, T_0/\varepsilon];(W^{2, \infty}_p(\Omega))^2) \\
& \cap C([0, T_0/\varepsilon];(W^{3, \infty}_p(\Omega))^2),
\end{split}
\end{equation*}
and
\begin{equation*}
(A)  \quad \left\|\frac{\partial^{r+s}}{\partial t^r\partial x^s}\Phi\right\|_{L^{\infty}([0, T_0/\varepsilon];(L^{\infty}(\Omega))^2)} \lesssim 1,\quad 0 \leq r \leq 3,\quad 0 \leq r+s \leq 3,\quad 0 < \varepsilon \leq 1,
\end{equation*}
where $W^{m, \infty}_p = \{u | u \in W^{m, \infty}(\Omega),\ \partial^l_x u(a)= \partial^l_x u(b), \ l=0, \cdots, m-1\}$ for $m \geq 1$, and here the boundary values are understood in the trace sense. In the subsequent discussion, we will properly omit $\Omega$ when referring to the space norm taken on $\Omega$. Denote $\Omega_{T_{\varepsilon}} = \Omega \times [0,T_0/\varepsilon]$, we assume the electromagnetic potentials $V \in C^2(\overline{\Omega}_{T_{\varepsilon}})$ and $A_1 \in C^2(\overline{\Omega}_{T_{\varepsilon}})$, and denote 
\begin{equation*}
(B) \quad V_{\max}:= \max_{(t, x) \in \overline{\Omega}_{T_{\varepsilon}}}|V(t, x)|,\quad A_{1, \max}:= \max_{(t, x) \in \overline{\Omega}_{T_{\varepsilon}}}|A_1(t, x)|.
\end{equation*}
Let $\Phi^n_j$ be the numerical approximations obtained from the FDTD methods and define the grid error function $\mathbf{e}^n = (\mathbf{e}^n_0, \mathbf{e}^n_1,\cdots,\mathbf{e}^n_N)^T \in X_N$ as:
\begin{equation}
\label{eq:e}
\mathbf{e}^n_j : = \Phi(t_n, x_j) - \Phi^n_j, \quad j=0,1, \cdots, N, \quad n\geq 0.
\end{equation}
We can establish the following error estimates for the FDTD methods under the corresponding stability conditions up to the time at $O(\varepsilon^{-1})$.

\begin{theorem}
\label{thm:FDTD}
Under the assumptions (A) and (B), there exist constants $h_0 > 0$ and $\tau_0 > 0$ sufficiently small and independent of $\varepsilon$, such that for any $0 < \varepsilon \leq 1$, when $0 < h \leq h_0$ and $0 < \tau < \tau_0$ and under the corresponding stability conditions, we have the following error estimate for the FDTD methods \eqref{eq:CNFD} - \eqref{eq:LFFD}
 with \eqref{eq:she0} and \eqref{eq:phi1}
\begin{equation}
\|\mathbf{e}^n\|_{l^2} \lesssim \frac{h^2}{\varepsilon} + \frac{\tau^2}{\varepsilon}, \quad 0 \leq n \leq \frac{T_0/\varepsilon}{\tau}.
\label{eq:FDTD_error}
\end{equation}
\end{theorem}

Based on the above theorem, the four FDTD methods share the same spatial/temporal resolution capacity for the Dirac equation \eqref{eq:Dirac_1D} up to the long time at $O(\varepsilon^{-1})$. In fact, given an accuracy bound $\delta > 0$, the $\varepsilon$-scalability of the FDTD methods should be taken as
\begin{equation}
h = O(\sqrt{\varepsilon\delta}) = O(\varepsilon^{1/2}),	 \quad \tau = O(\sqrt{\varepsilon\delta})  = O(\varepsilon^{1/2}),	\quad 0 < \varepsilon \leq 1.
\end{equation}
This implies that, in order to get ``correct'' numerical solution up to the time at $O(\varepsilon^{-1})$, one has to take the meshing strategy: $h=O(\varepsilon^{1/2})$ and  $\tau=O(\varepsilon^{1/2})$, which is very useful for practical computations on how to select mesh size and time step such that the numerical results are trustable.

For the physical observables including the total probability density and current density, weaker constraints are addmissible for the time-splitting method to solve the Schr\"odinger equation in the semiclassical regime \citep{BJP}. However, the convergence of the physical observables for the long-time dynamics of the Dirac equation are the same as that for the wave function as follows:
\begin{corollary}
Under the assumptions (A) and (B), with \eqref{eq:she0}, \eqref{eq:phi1} and corresponding stability conditions for CNFD, SIFD1, SIFD2 and LFFD, there exist constants $h_0 > 0$ and $\tau_0 > 0$ sufficiently small and independent of $\varepsilon$, such that for any $0 < \varepsilon \leq 1$, when $0 < h \leq h_0$ and $0 < \tau < \tau_0$, we have the following error estimates on the total probability density and current density for the FDTD methods \eqref{eq:CNFD} - \eqref{eq:LFFD}
\begin{equation}
\|\rho^n - \rho(t_n, \cdot)\|_{l^1} \lesssim \frac{h^2}{\varepsilon} + \frac{\tau^2}{\varepsilon},\quad \|\mathbf{J}^n - \mathbf{J}(t_n, \cdot)\|_{l^1} \lesssim \frac{h^2}{\varepsilon} + \frac{\tau^2}{\varepsilon}, \quad 0 \leq n \leq \frac{T_0/\varepsilon}{\tau},
\end{equation}
where $\rho^n$ and $\mathbf{J}^n$ are obtained from the wave function $\Phi^n$ through \eqref{eq:rho} and \eqref{eq:J} with $d=1$.
\end{corollary}

On the other hand, in order to obtain the error bounds for the FDFP methods in the long-time regime, we assume that there exists an integer $m_0 \geq 2$ such that the exact solution of the Dirac equation \eqref{eq:Dirac_1D} up to the time at $t = T_0/\varepsilon$ satisfies
\begin{equation*}
(C)\quad  \|\Phi\|_{L^{\infty}([0,T_0/\varepsilon];(H^{m_0}_p)^2)}	\lesssim 1, \quad \|\partial_t^s \Phi\|_{L^{\infty}([0,T_0/\varepsilon];(L^{\infty})^2)} \lesssim 1,\  s = 1, 2, 3,
\end{equation*}
where $H^k_p(\Omega) = \{u | u \in H^k(\Omega), \partial^l_x u(a) = \partial^l_x u(b), l=0, \cdots, k-1 \}$. In addition, we assume that the electromagnetic potentials satisfy
\begin{equation*}
(D)\quad  \|V\|_{W^{2,\infty}([0, T_0/\varepsilon]; L^{\infty})} +  \|A_1\|_{W^{2,\infty}([0, T_0/\varepsilon]; L^{\infty})} \lesssim 1.
\end{equation*}
Let $\Phi_j^n$ be the numerical approximations obtained from the FDFP methods, then we have the following theorem.

\begin{theorem}
\label{thm:FDFP}
Under the assumptions (C) and (D), and under the corresponding stability conditions, there exist constants $h_0, \tau_0 > 0$ sufficiently small and independent of $\varepsilon$ such that, for any $0 < \varepsilon \leq 1$, when $0 < h \leq h_0$ and $0 < \tau \leq \tau_0$, we have the following error estimate for the FDFP methods
\begin{equation}
\|\Phi(t_n, \cdot) - I_N(\Phi^n)\|_{L^2} \lesssim h^{m_0} + \frac{\tau^2}{\varepsilon}, \quad 0 \leq n \leq \frac{T_0/\varepsilon}{\tau}.	
\label{eq:FDFP_error}
\end{equation}
\end{theorem}

\subsection{Proof for FDTD}
Since the proofs for the error estimates of the FDTD methods are all similar, here we only show the details for CNFD. Other proofs can be derived in a similar manner.

\emph{Proof.}
Define the local truncation error $\widetilde{\xi}^n = (\widetilde{\xi}^n_0, \widetilde{\xi}^n_1, \cdots, \widetilde{\xi}^n_N)^T \in X_N$ of the CNFD method \eqref{eq:CNFD} with \eqref{eq:she0} for $0 \leq j \leq N$ and $n \geq 0$

\begin{align}
\label{eq:xi}	
\widetilde{\xi}^n_j := \ & i\delta^+_t\Phi(t_n, x_j) + i \sigma_1 \frac{\delta_x \Phi(t_n, x_j) + \delta_x \Phi(t_{n+1}, x_j)}{2} \nonumber\\
& - \left(\sigma_3 + \varepsilon V^{n+\frac{1}{2}}_j I_2 -\varepsilon A^{n+\frac{1}{2}}_{1, j} \sigma_1\right)\frac{\Phi(t_n, x_j) + \Phi(t_{n+1}, x_j)}{2}.
\end{align}
Applying in the Taylor expansion in \eqref{eq:xi}, we obtain for $j = 0, 1, \cdots, N-1$ and $n \geq 0$,
\begin{equation*}
\begin{split}
\widetilde{\xi}^n_j = \ &\frac{i\tau^2}{6} \partial_{ttt} \Phi(t'_n, x_j) + \frac{ih^2}{12}\sigma_1\partial_{xxx}\Phi(t_n, x'_j) + \frac{ih^2}{12}\sigma_1\partial_{xxx}\Phi(t_{n+1}, x''_j) + \frac{i\tau^2}{4}\sigma_1\partial_{xtt}\Phi(t''_n, x_j)\\
&  -\frac{\varepsilon \tau^2}{8}\Big(\partial_{tt}V(t_n, x_j) I_2 - \partial_{tt}A_{1}(t_n, x_j) \sigma_1\Big) \Phi(t'''_n, x_j) - \frac{\tau^2}{4}\Big(\sigma_3 + \varepsilon V^{n+\frac{1}{2}}_j I_2 -\varepsilon A^{n+\frac{1}{2}}_{1, j} \sigma_1\Big)\partial_{tt}\Phi(t''''_n, x_j),  
\end{split}
\end{equation*}
where $t_{n-1} < t'_n, t''_n, t'''_n, t''''_n < t_{n+1}$, $x_{j-1} < x'_j, x''_j < x_{j+1}$. Noticing \eqref{eq:Dirac_1D} and the assumptions $(A)$ and $(B)$, we can get the following estimates for $n \geq 1$
\begin{equation}
\label{eq:xib}
\|\widetilde{\xi}^n\|_{l^{\infty}}	= \max_{0 \leq j \leq N-1} |\xi^n_j| \lesssim h^2 + \tau^2, \quad \|\widetilde{\xi}^n\|_{l^2} \lesssim \|\widetilde{\xi}^n\|_{l^{\infty}} \lesssim h^2 + \tau^2.
\end{equation}
By the definition of the error function \eqref{eq:e}, we obtain for $0 \leq j \leq N-1$ and $n \geq 0$, 
\begin{equation}
i\delta^+_t\mathbf{e}^n_j = -i\sigma_1\delta_x \mathbf{e}^{n+1/2}_j + \left(\sigma_3 + \varepsilon V^{n+\frac{1}{2}}_j I_2 -\varepsilon A^{n+\frac{1}{2}}_{1, j} \sigma_1\right)\mathbf{e}^{n+1/2}_j + \widetilde{\xi}^n_j,
\end{equation}
where the initial and boundary conditions are given as 
\begin{equation}
\label{eq:error}
\mathbf{e}^n_0 = \mathbf{e}^n_N,\quad \mathbf{e}^n_{-1} = \mathbf{e}^n_{N-1},\quad n \geq 0, \quad \mathbf{e}^0_j = 0, \quad j = 0, 1, \cdots, N.
\end{equation}
Multiplying \eqref{eq:error} by the left by $\tau h \left(\mathbf{e}^{n+1}_j + \mathbf{e}^n_j\right)^{\ast}$, taking the imaginary part, then summing up for $j = 0, 1, \cdots, N-1$, using Cauchy inequality, noticing \eqref{eq:xib}, we can get
\begin{align}
\|\mathbf{e}^{n+1}\|^2_{l^2} - \|\mathbf{e}^{n}\|^2_{l^2} &= \tau h \text{Im}\left[\sum^{N-1}_{j = 0}\left(\mathbf{e}^{n+1}_j + \mathbf{e}^n_j\right)^{\ast}\widetilde{\xi}^n_j\right] \nonumber\\
&\lesssim \varepsilon \tau \left(\|\mathbf{e}^{n+1}\|^2_{l^2} + \|\mathbf{e}^{n}\|^2_{l^2}\right) + \frac{\tau}{\varepsilon} \left(h^2 + \tau^2\right)^2.
\end{align}
Then summing up the above inequality from $n = 0$ to $n = m-1$ gives
\begin{equation}
\|{\mathbf{e}^{m}}\|^2_{l^2} \lesssim\eps\tau\sum_{k=0}^{m}\|\mathbf{e}^k\|^2_{l^2} + \frac{m\tau}{\eps}(h^2 + \tau^2)^2.
\end{equation}
For small enough $\tau_0$, taking into account the fact that $m\leq \frac{T_0/\varepsilon}{\tau}$, we obtain
\begin{equation}
\|\mathbf{e}^m\|^2_{l^2} \leq T_0e^{T_0}\left(\frac{h^2}{\varepsilon} + \frac{\tau^2}{\varepsilon}\right)^2,
\end{equation}
by applying the discrete Gronwall's inequality, which implies the error bound \eqref{eq:FDTD_error}.
\hfill $\square$ \bigskip

\subsection{Proof for FDFP}
We only show the proof of the long-time error bound for the LFFP method. Proofs for the other three FDFP methods are similar and we omit the details here for brevity. 

\emph{Proof.}
For $0 \leq n \leq \frac{T_0/\varepsilon}{\tau}$, we rewrite the error as 
\begin{equation}
\Phi(t_n, \cdot) - I_N(\Phi^n)	= \Phi(t_n, \cdot)  - P_N(\Phi(t_n, \cdot)) + P_N(\Phi(t_n, \cdot)) - I_N(\Phi^n)	
\end{equation}
The regularity of $\Phi(t_n, \cdot)$ implies
\begin{equation}
\|\Phi(t_n, \cdot)  - P_N(\Phi(t_n, \cdot))\|_{L^2} \lesssim h^{m_0}.
\end{equation}
Thus, it remains to establish the error bound for the error
\begin{equation}
\eta^n := P_N(\Phi(t_n, \cdot)) - I_N(\Phi^n), \quad	0 \leq n \leq \frac{T_0/\varepsilon}{\tau}.
\end{equation}

For $n=0$, from the discretization of the initial data, noticing \eqref{eq:LFFS_initial} and \eqref{eq:e}, we have 
\begin{equation}
\|\eta^0\|_{L^2}  = \|P_N(\Phi_0) - I_N(\Phi^0)\|_{L^2} \lesssim h^{m_0},
\end{equation}
which implies the error bound \eqref{eq:FDFP_error} is valid for $n = 0$.

Taking the Fourier transform of the LFFP method \eqref{eq:LFFS}, for $l \in \mathcal{T}_N$, we obtain 
\begin{equation}
\label{eq:In}
i \frac{(\widetilde{\Phi^{n+1}})_l - (\widetilde{\Phi^{n-1}})_l}{2\tau} = \sigma_1\mu_l  (\widetilde{\Phi^{n}})_l + \sigma_3  (\widetilde{\Phi^{n}})_l + \varepsilon \left( (\widetilde{V^n\Phi^{n}})_l - \sigma_1(\widetilde{A^n_{1}\Phi^{n}})_l \right).
\end{equation}
Define the local truncation error $\xi^n \in Y_N$ as
\begin{equation}
\xi^n = \sum_{l \in \mathcal{T}_N} (\widehat{\xi^n})_l e^{i\mu_l(x-a)},	
\end{equation}
with
\begin{equation}
\begin{aligned}
(\widehat{\xi^0})_l & := i\delta_t^+ (\widehat{\Phi_0})_l {- \sigma_1} \mu_l  (\widehat{\Phi_0})_l -\sigma_3(\widehat{\Phi_0})_l -\varepsilon \left( (\widehat{V^0\Phi_0})_l- \sigma_1(\widehat{A^0_1\Phi_0})_l \right),  \\	
(\widehat{\xi^n})_l & := i\delta_t \widehat{\Phi}_l (t_n) { -\sigma_1} \mu_l  \widehat{\Phi}_l (t_n) -\sigma_3\widehat{\Phi}_l (t_n) -\varepsilon \left( (\widehat{V^n\Phi(t_n)})_l- \sigma_1(\widehat{A^n_1\Phi(t_n)})_l \right), \ n \geq 1.
\label{eq:F1}
\end{aligned}
\end{equation}

Similar to the proof of CNFD, applying the Taylor expansion in \eqref{eq:F1}, under the assumptions (C) and (D), we have
\begin{equation}
\left|(\widehat{\xi^0})_l\right| \lesssim  \tau, \quad \left|(\widehat{\xi^n})_l\right| \lesssim  \tau^2, \quad l\in \mathcal{T}_N,\quad n\geq 1,
\end{equation}
which immediately implies $\|\xi^0\|_{L^2}\leq \tau$, and for $n \geq 1$
\begin{equation}
\|\xi^n\|_{L^2} \lesssim  \tau^2.
\label{eq:xi_error}
\end{equation}

We denote the error of the electromagnetic potentials
\begin{equation}
R^n_j = \sum_{l \in \mathcal{T}_N}(\widehat{R^n})_l e^{2ijl\pi/N},	
\end{equation}
with
\begin{equation}
(\widehat{R^n})_l = \varepsilon \left[(\widehat{V^n\Phi(t_n)})_l- \sigma_1(\widehat{A^n_1\Phi(t_n)})_l \right]-\varepsilon\left[(\widehat{I_N(V^n\Phi^n}))_l -\sigma_1(\widehat{I_N(A^n_1\Phi^n}))_l\right],\quad l \in \mathcal{T}_N,\quad n \geq 0.
\end{equation}
For each $l \in \mathcal{T}_N$, subtracting \eqref{eq:In} from \eqref{eq:F1}, we get
\begin{align}
\label{eq:LFFS_ef}
& i\delta_t (\widehat{\eta^n})_l = \sigma_1 \mu_l (\widehat{\eta^n})_l + \sigma_3 (\widehat{\eta^n})_l  + (\widehat{\xi^n})_l + (\widehat{R^n})_l, \quad 1 \leq n \leq \frac{T_0/\varepsilon}{\tau}-1, \nonumber \\
& (\widehat{\eta^0})_l = 0, \quad (\widehat{\eta^1})_l = -i\tau (\widehat{\xi^0})_l, \quad l \in \mathcal{T}_N.
\end{align}
The error equation \eqref{eq:LFFS_ef} and the estimate $\|\xi^0\|_{L^2}\leq \tau$ imply
\begin{equation}
\|\eta^1\|_{L^2} = \tau \|\xi^0\|_{L^2} \lesssim  \tau^2.	
\end{equation}
Thus, we immediately obtain
\begin{equation}
\|\Phi(t_1, \cdot) - I_N(\Phi^1)\|_{L^2} \lesssim h^{m_0} + \tau^2,
\end{equation}
which indicates that the error bound \eqref{eq:FDFP_error} is valid when $n=1$.

Now, we assume that \eqref{eq:FDFP_error} is valid for all $1 \leq n \leq m-1 \leq \frac{T_0/\varepsilon}{\tau}-1$, then we need to show that it is still valid when $n = m$.
Under the assumptions (C) and (D), we have
\begin{equation}
\|R^n\|_{L^2} \lesssim \varepsilon h^{m_0} + \tau^2, \quad 1 \leq n \leq m-1.	
\end{equation}

Define the energy $\mathcal{E}^n$( $n = 0, 1, \cdots$)  as
\begin{align}
\mathcal{E}^{n+1}= & \; \|{\eta}^{n+1}\|^2_{L^2}+\|{\eta}^{n}\|^2_{L^2}	-2\tau\text{Im}\left(h\sum_{l\in\mathcal{T}_N}\mu_l(\widehat{{\eta}^{n+1}})_l^{\ast}\sigma_1 (\widehat{\eta^n})_l \right).
\end{align}
Under the stability condition \eqref{eq:st_LFFS} $\tau \leq \frac{\tau_1h}{2(V_{\max}h+\sqrt{h^2+(\pi+ h A_{1, \max})^2})}$, if we take $\tau_1 = \frac{1}{4}$, we could derive  $\tau (1+\varepsilon(V_{\max} + A_{1, \max})) \leq \frac{1}{8}$, which gives
\begin{equation}
\frac{1}{2}\left(\|{\eta}^{n+1}\|^2_{L^2}+\|{\eta}^{n}\|^2_{L^2}\right) \leq \mathcal{E}^{n+1} \leq \frac{3}{2}\left(\|{\eta}^{n+1}\|^2_{L^2}+\|{\eta}^{n}\|^2_{L^2}\right), \quad n \geq 0,
\label{eq:eE}
\end{equation}
by using Cauchy inequality.

Furthermore, multiply $2\tau h\left((\widehat{\eta^{n+1}})_l +(\widehat{\eta^{n-1}})_l\right)^{\ast}$ from left on both sides to \eqref{eq:LFFS_ef}, then sum up for $l \in \mathcal{T}_N$, and take the imaginary part, similar to \eqref{al:massproof}, we have
\begin{align}
\mathcal{E}^{n+1} - \mathcal{E}^n  = & 2 \tau h \text{Im} \left[ \sum_{l\in\mathcal{T}_N}\left( (\widehat{\eta^{n+1}})_l + (\widehat{\eta^{n-1}})_l\right)^{\ast}\left((\widehat{\xi^n})_l+(\widehat{R^n})_l\right)\right] \nonumber \\
\lesssim &\;  \varepsilon \tau \left(\|{\eta}^{n+1}\|^2_{L^2} + \|{\eta}^{n-1}\|^2_{L^2}\right) + \frac{\tau}{\varepsilon} \left(\|\xi^n\|^2_{L^2}+ \|R^n\|^2_{L^2}\right)\nonumber \\
\lesssim &\; \varepsilon \tau \left(\mathcal{E}^{n+1}+\mathcal{E}^{n}\right) + \frac{\tau}{\varepsilon}\left( \varepsilon h^{m_0}+\tau^2\right)^2,
\end{align}
by noticing \eqref{eq:xi_error}. Summing the above inequalities from $1$ to  $m-1$, we get 
\begin{equation}
\mathcal{E}^m-\mathcal{E}^1	\lesssim \varepsilon \tau \sum^m_{k=1}\mathcal{E}^k + \frac{m\tau}{\varepsilon} \left(\varepsilon h^{m_0}+\tau^2\right)^2,\quad 1 \leq m \leq \frac{T_0/\varepsilon}{\tau}.
\label{eq:E_m}
\end{equation}
So if we take $\tau_0$ sufficiently small, under the discrete Gronwall's inequality for \eqref{eq:E_m}, it can be obtained that
\begin{equation}
\mathcal{E}^m \lesssim \left(h^{m_0}+ \frac{\tau^2}{\varepsilon}\right)^2,
\end{equation}
which immediately implies 
\begin{equation}
\|\eta^m\|_{L^2} \lesssim h^{m_0}+ \frac{\tau^2}{\varepsilon}.
\end{equation}
Thus, we have 
\begin{equation}
\|\Phi(t_m, \cdot) - I_N(\Phi^m)\|_{l^2} \lesssim h^{m_0} + \frac{\tau^2}{\varepsilon},
\end{equation}
which indicates the error bound \eqref{eq:FDFP_error} is still valid when $n = m$. Hence, the proof of Theorem \ref{thm:FDFP} by the method of mathematical induction.
\hfill $\square$ \bigskip

\begin{remark}
For the FDTD methods, the accuracy and $\varepsilon$-scalability in space could be improved by using higher order finite difference discretizations, but the spatial errors will still depend on the small parameter $\varepsilon$ in the long-time regime. Similarly, the spatial error bounds of the finite element and finite volume discretizations in space also depend on the parameter $\varepsilon$ in the long-time regime. Uniform error bounds in space for the long-time dynamics can be achieved by the spectral method. For the spectral method, the computation of  spatial derivative operator is carried out in the Fourier space. In contrast to the finite difference, finite element and finite volume methods, the spectral approximations of the differential operator are exact for all the Fourier modes \citep{TE}. From the proof of the error bounds for the FDTD methods, we can observe that the errors of the approximations for the spatial differential operator accumulate in the long-time regime, which results in the fact the spatial errors depend on the small parameter $\varepsilon$. By using the spectral discretization in space, the spatial errors only come from the projection and the potential term. In the long-time regime, the accumulation of the spatial errors from the potential is independent of the small parameter $\varepsilon$ since the strength of the potential is at $O(\varepsilon)$. In summary, the spectral method performs much better than other spatial discretizations for the long-time dynamics of the Dirac equation with small potential.
\end{remark}

\section{Numerical results}
In this section, we compare the spatial and temporal resolutions of different numerical methods including the FDTD methods and the FDFP methods of the Dirac equation in 1D. We pay particular attention to the $\varepsilon$-scalability of different methods in the long-time regime, i.e., $0 < \varepsilon \ll 1$. 

\subsection{Spatial/Temporal resolution for the periodic problem}
In this numerical experiment, the problem is solved numerically on a torus $\Omega = (0, 2\pi)$ with periodic boundary conditions. We choose the electromagnetic potentials in the Dirac equation \eqref{eq:Dirac_21} as
\begin{equation}
V(t, x) = \frac{1}{1+\sin^2(x)},\quad A_1(x) = \cos(x) + \sin(2x), \quad x \in \Omega,	
\end{equation}
and the initial data are taken as
\begin{equation}
\phi_1(0, x) = \frac{1}{1+\sin^2(x)},\quad \phi_2(0, x) = \frac{1}{3+\cos(x)},\quad x \in \Omega.
\end{equation}
Since the exact solution of the Dirac equation is unknown, we use the time splitting Fourier pseudospectral (TSFP) method \citep{BCJT} with a fine mesh size $h_e = \pi/1024$  and a very small time step $\tau_e = 10^{-4}$ to get the `reference' solution numerically for comparison of the results from FDTD and FDFP methods. In order to quantify the numerical errors, we introduce the following discrete $l^2$ errors of the wave function $\Phi$, the total probability density $\rho$ and the current density $\mathbf{J}$
\begin{equation*}
e_{\Phi}(t_n) = \|\Phi^n - \Phi(t_n, \cdot)\|_{l^2}, \quad e_{\rho}(t_n) = \|\rho^n - \rho(t_n, \cdot)\|_{l^1}, \quad e_{\mathbf{J}}(t_n) = \|{\mathbf{J}}^n - \mathbf{J}(t_n, \cdot)\|_{l^1},
\end{equation*}
where $\rho^n$ and $\mathbf{J}^n$ can be computed from the numerical solution $\Phi$ via \eqref{eq:rho} and \eqref{eq:J} with $d=1$, respectively.

\begin{table}[h!]
\caption{Spatial and temporal errors of the CNFD method for the wave function $e_{\Phi}(t=2/\varepsilon)$ of the Dirac equation \eqref{eq:Dirac_21} in 1D.}
\label{tab:CNFD_phi}
\renewcommand{\arraystretch}{1.3}
\centering
\vspace{0.3cm}
\begin{tabular}{cccccc}
\hline
$e_{\Phi}(t=2/\varepsilon)$ &$\begin{matrix} h_0 = \frac{\pi}{64} \\ \tau_0 = 0.05 \end{matrix}$ & $\begin{matrix} h_0/2 \\ \tau_0/2 \end{matrix}$ & $\begin{matrix} h_0/2^2 \\ \tau_0/2^2 \end{matrix}$ &  $\begin{matrix} h_0/2^3 \\ \tau_0/2^3 \end{matrix}$ &  $\begin{matrix} h_0/2^4 \\ \tau_0/2^4 \end{matrix}$   \\
\hline
$\varepsilon = 1$ & \bf{3.40E-2} & 8.53E-3 & 2.14E-3 & 5.34E-4 & 1.33E-4  \\
order & \bf{-} & 1.99 & 1.99 & 2.00 & 2.01  \\
\hline
$\varepsilon = 1/4$ & 4.02E-2 & \bf{1.03E-2} & 2.57E-3 & 6.43E-4 & 1.61E-4  \\
order & - & \bf{1.96} & 2.00 & 2.00 & 2.00  \\
\hline
$\varepsilon = 1/4^2$ & 1.19E-1 & 3.34E-2 & \bf{8.50E-3} & 2.13E-3 & 5.33E-4  \\
order & - & 1.83 & \bf{1.97} & 2.00 & 2.00 \\
\hline
$\varepsilon = 1/4^3$ & 3.04E-1 & 1.17E-1 & 3.27E-2 & \bf{8.31E-3} & 2.08E-3  \\
order & - & 1.38 & 1.84 & \bf{1.98} & 2.00 \\
\hline
$\varepsilon = 1/4^4$ & 8.75E-1 & 3.02E-1 & 1.16E-1 & 3.26E-2 & \bf{8.27E-3}  \\
order & - & 1.53 & 1.38 & 1.83 & \bf{1.98} \\
\hline
\end{tabular}
\end{table}

\begin{table}[h!]
\caption{Spatial and temporal errors of the CNFD method for the total probability $e_{\rho}(t=2/\varepsilon)$ of the Dirac equation \eqref{eq:Dirac_21} in 1D.}
\label{tab:CNFD_rho}
\renewcommand{\arraystretch}{1.3}
\centering
\vspace{0.3cm}
\begin{tabular}{cccccc}
\hline
$e_{\rho}(t=2/\varepsilon)$ &$\begin{matrix} h_0 = \frac{\pi}{64} \\ \tau_0 = 0.05 \end{matrix}$ & $\begin{matrix} h_0/2 \\ \tau_0/2 \end{matrix}$ & $\begin{matrix} h_0/2^2 \\ \tau_0/2^2 \end{matrix}$ &  $\begin{matrix} h_0/2^3 \\ \tau_0/2^3 \end{matrix}$ &  $\begin{matrix} h_0/2^4 \\ \tau_0/2^4 \end{matrix}$   \\
\hline
$\varepsilon = 1$  & \bf{3.67E-2} & 8.87E-3 & 2.20E-3 & 5.48E-4 & 1.37E-4  \\
order & \bf{-} & 2.05 & 2.01 & 2.01 & 2.00  \\
\hline
$\varepsilon = 1/4$  & 6.38E-2 & \bf{1.54E-2} & 3.75E-3 & 9.30E-4 & 2.32E-4  \\
order & - & \bf{2.05} & 2.04 & 2.01 & 2.00  \\
\hline
$\varepsilon = 1/4^2$  & 2.60E-1 & 6.30E-2 & \bf{1.54E-2} & 3.81E-3 & 9.50E-4  \\
order & - & 2.05 & \bf{2.03} & 2.02 & 2.00 \\
\hline
$\varepsilon = 1/4^3$  & 4.80E-1 & 1.93E-1 & 5.03E-2 & \bf{1.21E-2} & 3.04E-3  \\
order & - & 1.31 & 1.94 & \bf{2.06} & 1.99 \\
\hline
$\varepsilon = 1/4^4$  & 1.41 & 8.34E-1 & 3.21E-1 & 8.05E-2 & \bf{2.08E-2}  \\
order & - & 0.76 & 1.38 & 2.00 & \bf{1.95} \\
\hline
\end{tabular}
\end{table}

\begin{table}[h!]
\caption{Spatial and temporal errors of the CNFD method for the current density $e_{\mathbf{J}}(t=2/\varepsilon)$ of the Dirac equation \eqref{eq:Dirac_21} in 1D.}
\label{tab:CNFD_J}
\renewcommand{\arraystretch}{1.3}
\centering
\vspace{0.3cm}
\begin{tabular}{cccccc}
\hline
$e_{\mathbf{J}}(t=2/\varepsilon)$ &$\begin{matrix} h_0 = \frac{\pi}{64} \\ \tau_0 = 0.05 \end{matrix}$ & $\begin{matrix} h_0/2 \\ \tau_0/2 \end{matrix}$ & $\begin{matrix} h_0/2^2 \\ \tau_0/2^2 \end{matrix}$ &  $\begin{matrix} h_0/2^3 \\ \tau_0/2^3 \end{matrix}$ &  $\begin{matrix} h_0/2^4 \\ \tau_0/2^4 \end{matrix}$   \\
\hline
$\varepsilon = 1$  & \bf{9.75E-2} & 2.44E-2 & 6.09E-3 & 1.52E-3 & 3.81E-4  \\
order & \bf{-} & 2.00 & 2.00 & 2.00 & 2.00  \\
\hline
$\varepsilon = 1/4$  & 9.91E-2 & \bf{2.43E-2} & 6.10E-3 & 1.53E-3 & 3.82E-4  \\
order & - & \bf{2.03} & 1.99 & 2.00 & 2.00  \\
\hline
$\varepsilon = 1/4^2$  & 4.05E-1 & 1.30E-1 & \bf{3.42E-2} & 8.59E-3 & 2.15E-3  \\
order & - & 1.64 & \bf{1.93} & 1.99 & 2.00 \\
\hline
$\varepsilon = 1/4^3$  & 1.32 & 3.87E-1 & 1.32E-1 & \bf{3.47E-2} & 8.69E-3  \\
order & - & 1.77 & 1.55 & \bf{1.93} & 2.00 \\
\hline
$\varepsilon = 1/4^4$  & 3.39 & 9.31E-1 & 3.38E-1 & 1.06E-1 & \bf{2.73E-2}  \\
order & - & 1.86 & 1.46 & 1.67 & \bf{1.96} \\
\hline
\end{tabular}
\end{table}

\begin{table}[h!]
\caption{Spatial errors of the LFFP method for the wave function $e_{\Phi}(t=2/\varepsilon)$ of the Dirac equation \eqref{eq:Dirac_21} in 1D.}
\label{tab:LFFP_h}
\renewcommand{\arraystretch}{1.3}
\centering
\vspace{0.3cm}
\begin{tabular}{ccccc}
\hline
$e_{\Phi}(t=2/\varepsilon)$ &$h_0 = \frac{\pi}{4}$ & $h_0/2 $ &$h_0/2^2 $ & $h_0/2^3$  \\
\hline
$\varepsilon = 1$  & 3.97E-1 & 1.43E-2 & 6.35E-6 & 1.04E-9 \\
\hline
$\varepsilon = 1/2$  & 2.13E-1 & 7.21E-3 & 6.31E-6 & 1.16E-9 \\
\hline
$\varepsilon = 1/2^2$  & 1.26E-1 & 3.08E-3 & 1.70E-6 & 1.46E-9 \\
\hline
$\varepsilon = 1/2^3$  & 4.33E-2 & 1.14E-3 & 2.42E-6 & 2.36E-9 \\
\hline
$\varepsilon = 1/2^4$  & 3.70E-2 & 1.40E-3 & 2.31E-6 & 4.50E-9 \\
\hline
\end{tabular}
\end{table}

\begin{table}[h!]
\caption{Temporal errors of the LFFP method for the wave function $e_{\Phi}(t=2/\varepsilon)$ of the Dirac equation \eqref{eq:Dirac_21} in 1D.}
\label{tab:LFFP_t}
\renewcommand{\arraystretch}{1.3}
\centering
\vspace{0.3cm}
\begin{tabular}{cccccc}
\hline	
$e_{\Phi}(t=2/\varepsilon)$ &$\tau_0 = 0.01$ & $\tau_0/2 $ &$\tau_0/2^2 $ & $\tau_0/2^3$ & $\tau_0/2^4$ \\
\hline
$\varepsilon = 1/2^6$  & \bf{1.72E-2} & 4.30E-3 & 1.08E-3 & 2.69E-4 & 6.72E-5  \\
order & \bf{- }& 2.00 & 1.99 & 2.01 & 2.00  \\
\hline
$\varepsilon = 1/2^8$  & 6.61E-2 & \bf{1.70E-3} & 4.27E-3 & 1.07E-3 & 2.67E-4  \\
order & - & \bf{1.96} & 1.99 & 2.00 & 2.00  \\
\hline
$\varepsilon = 1/2^{10}$  & 2.26E-1 & 6.60E-2 & \bf{1.70E-2} & 4.27E-3 & 1.07E-3 \\
order & - & 1.78 & \bf{1.96} & 1.99 & 2.00 \\
\hline
$\varepsilon = 1/2^{12}$  & 6.77E-1 & 2.26E-1 & 6.59E-2 & \bf{1.70E-2} & 4.26E-3 \\
order & - & 1.58 & 1.78 & \bf{1.95} & 2.00 \\
\hline
$\varepsilon = 1/2^{14}$  & 1.11 & 6.77E-1 & 2.26E-1 & 6.59E-2 & \bf{1.70E-2}  \\
order & - & 0.71 & 1.58 & 1.78 & \bf{1.95}\\
\hline
\end{tabular}
\end{table}

From Tables \ref{tab:CNFD_phi}-\ref{tab:CNFD_J} and additional similar numerical results not shown here for brevity, we can draw the following observations:

(i) For any fixed $\varepsilon = \varepsilon_0 > 0$, the FDTD methods are uniformly second-order accurate in both space and time, which agree with those results in the literature.

(ii) In the long-time regime, the second order convergence in space and time of the FDTD methods can only be observed when $0 < h \lesssim \varepsilon^{1/2}$ and $0 < \tau \lesssim \varepsilon^{1/2}$ (cf. upper triangles above the diagonals in Tables \ref{tab:CNFD_phi} -\ref{tab:CNFD_J}, where the diagonal lines in bold correspond to $h \sim \varepsilon^{1/2}$ and $\tau \sim \varepsilon^{1/2}$), which confirm our error estimates in Theorem \ref{thm:FDTD} and show that they are sharp.

From Tables \ref{tab:LFFP_h} - \ref{tab:LFFP_t} and additional numerical results for the other three FDFP methods, we can draw the following observations on the FDFP methods for the Dirac equation up to the time at $O(1/\varepsilon)$:

(i) In space, the FDFP methods are uniformly spectral accurate for any $0 < \varepsilon \leq 1$ (cf. each row in Table \ref{tab:LFFP_h}) and the spatial errors are almost independent of $\varepsilon$ (cf. each column in Table \ref{tab:LFFP_h}).

(ii) In time, for any fixed $\varepsilon = \varepsilon_0 > 0$, the FDFP methods show second order convergence in time (cf. first row in Table \ref{tab:LFFP_t}), which agree with those results in the literature. While, in the long-time time regime up to the time at $O(1/\varepsilon)$, the second order convergence in time can only be observed when $0 < \tau \lesssim \varepsilon^{1/2}$ (cf. the upper triangle above the diagonal in Table \ref{tab:LFFP_t}, where the bold diagonal line corresponds to $\tau \sim \varepsilon^{1/2}$), which confirm our error estimates in Theorem \ref{thm:FDFP}, and show that they are sharp.

\begin{figure}[ht!]
\begin{minipage}{0.5\textwidth}
\centerline{\includegraphics[width=8cm,height=6cm]{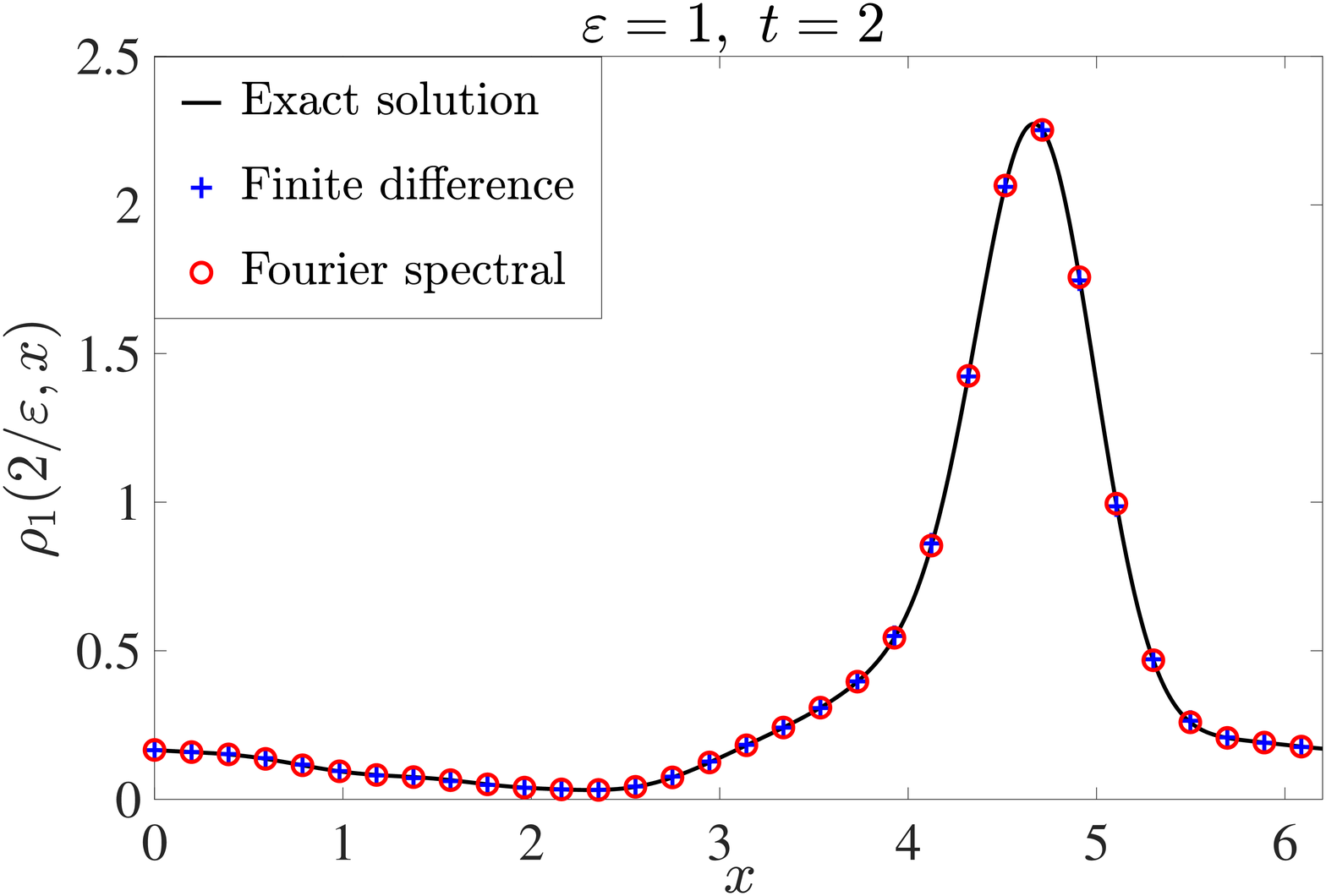}}
\end{minipage}
\begin{minipage}{0.5\textwidth}
\centerline{\includegraphics[width=8cm,height=6cm]{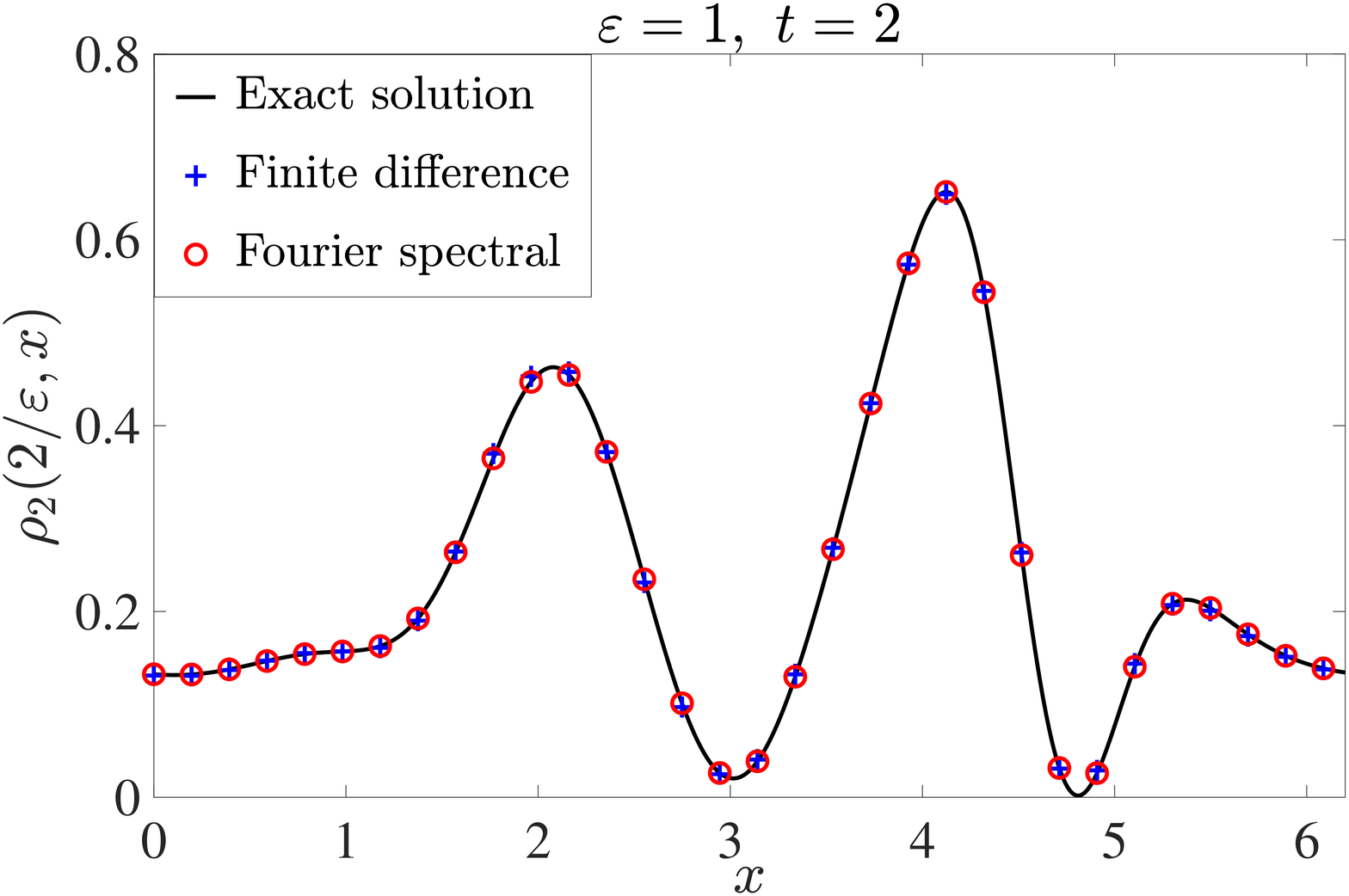}}
\end{minipage}
\begin{minipage}{0.5\textwidth}
\centerline{\includegraphics[width=8cm,height=6cm]{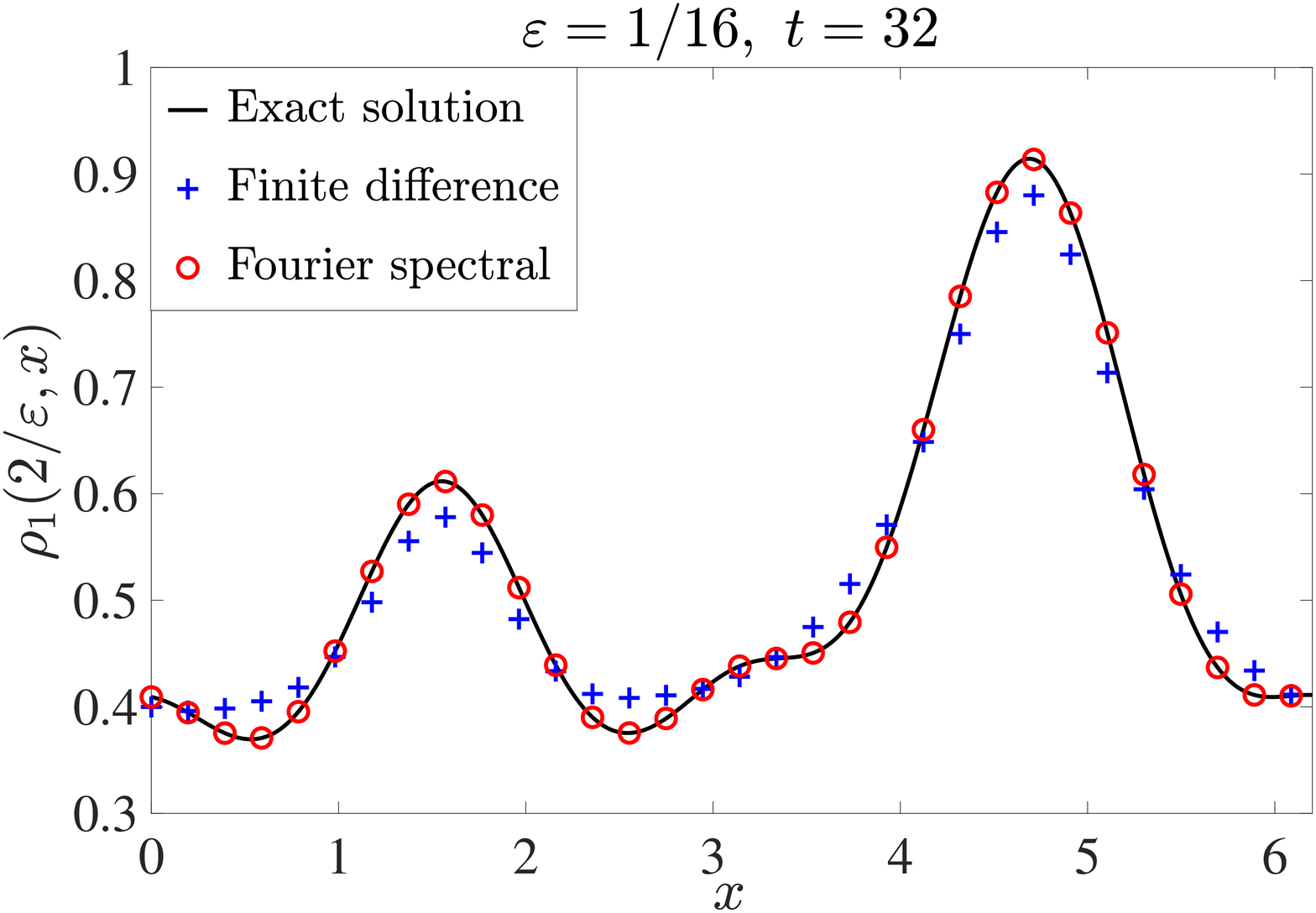}}
\end{minipage}
\begin{minipage}{0.5\textwidth}
\centerline{\includegraphics[width=8cm,height=6cm]{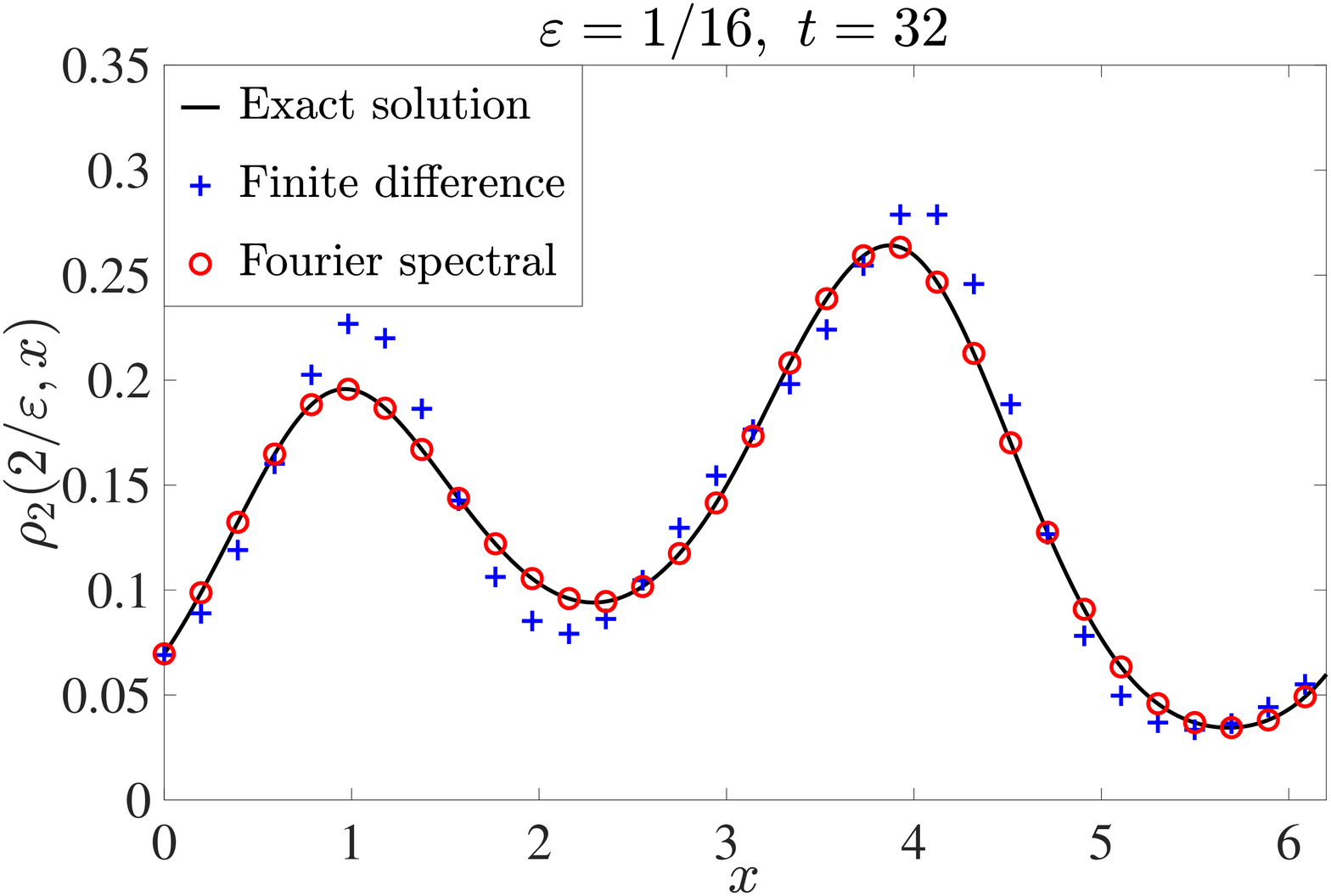}}
\end{minipage}
\begin{minipage}{0.5\textwidth}
\centerline{\includegraphics[width=8cm,height=6cm]{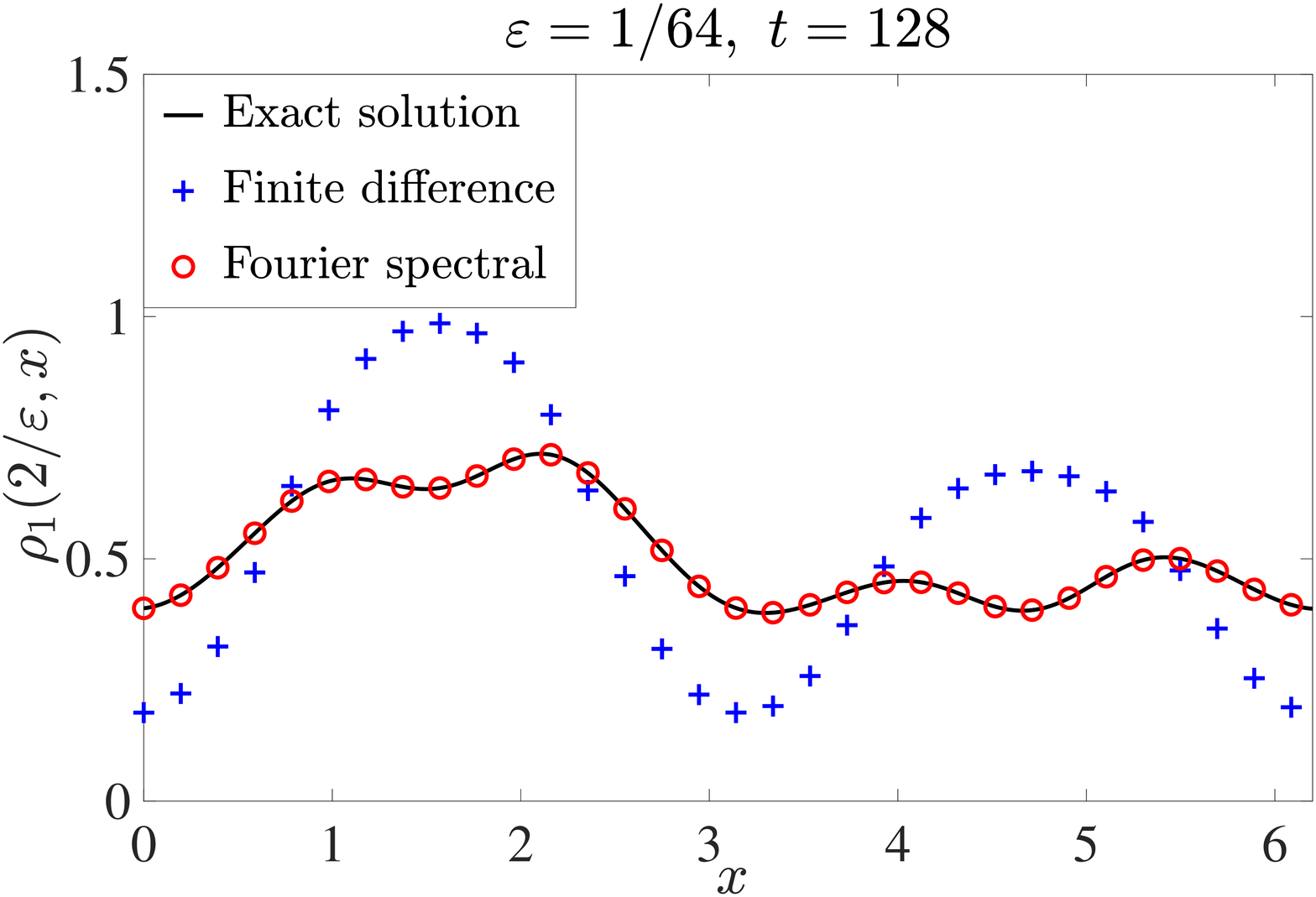}}
\end{minipage}
\begin{minipage}{0.5\textwidth}
\centerline{\includegraphics[width=8cm,height=6cm]{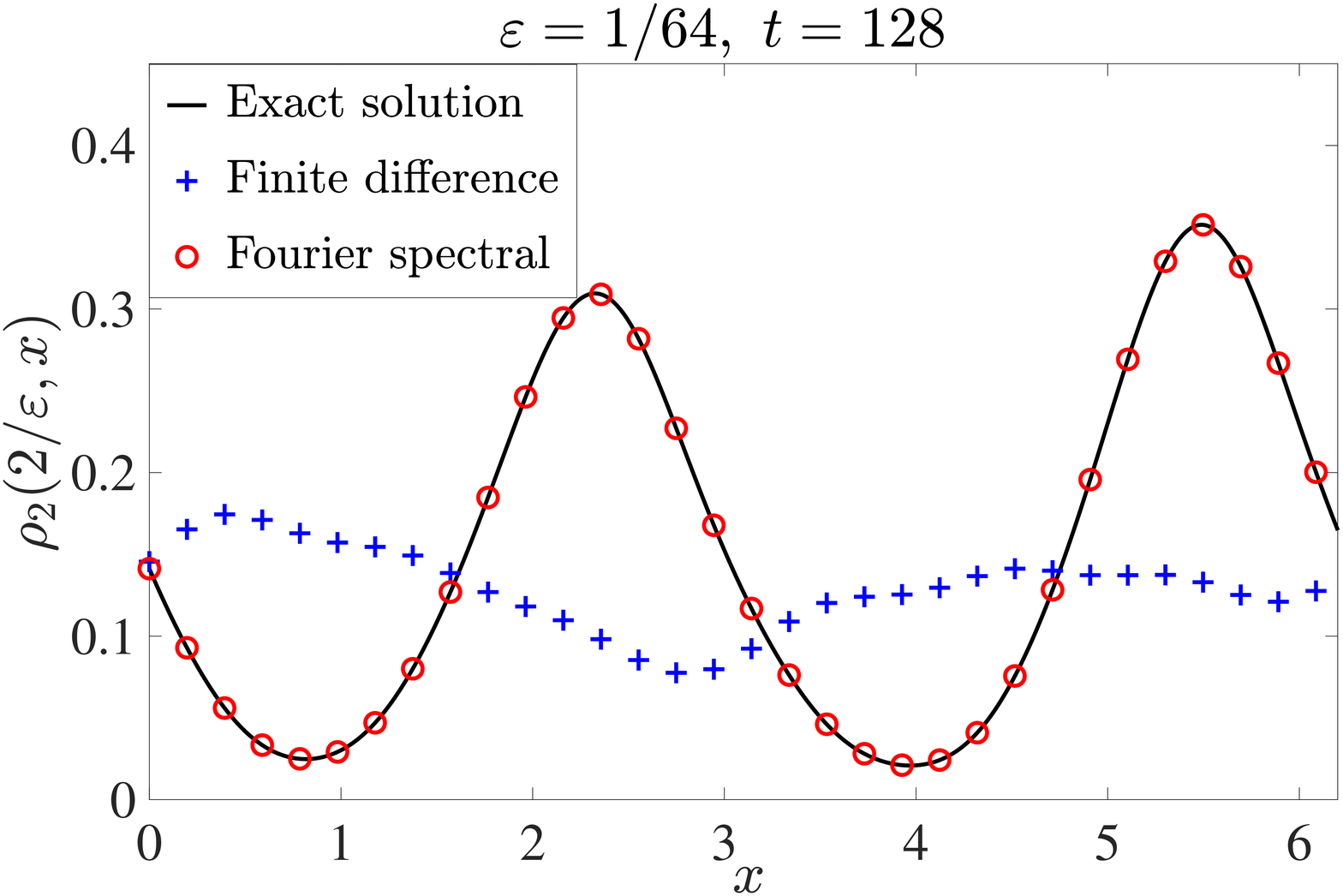}}
\end{minipage}
\caption{Long-time dynamics of the densities $\rho_1(2/\varepsilon, x) = |\phi_1(2/\varepsilon, x)|^2$ (left) and $\rho_2(2/\varepsilon, x) = |\phi_2(2/\varepsilon, x)|^2$ (right) of the Dirac equation  \eqref{eq:Dirac_21} in 1D with different $\varepsilon$.}
\label{fig:comp}
\end{figure}

\begin{figure}[ht!]
\begin{minipage}{0.5\textwidth}
\centerline{\includegraphics[width=8cm,height=6cm]{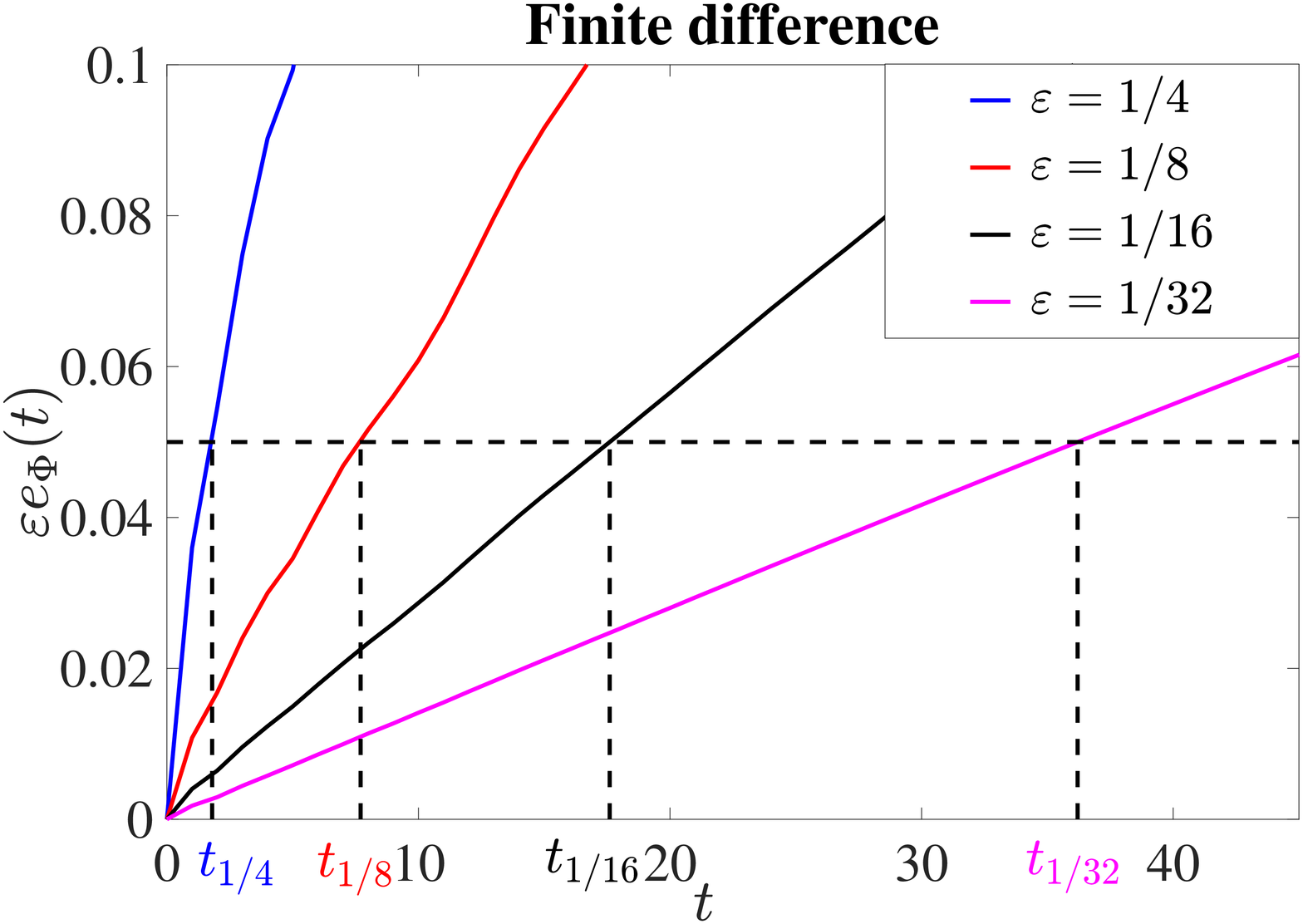}}
\end{minipage}
\begin{minipage}{0.5\textwidth}
\centerline{\includegraphics[width=8cm,height=6cm]{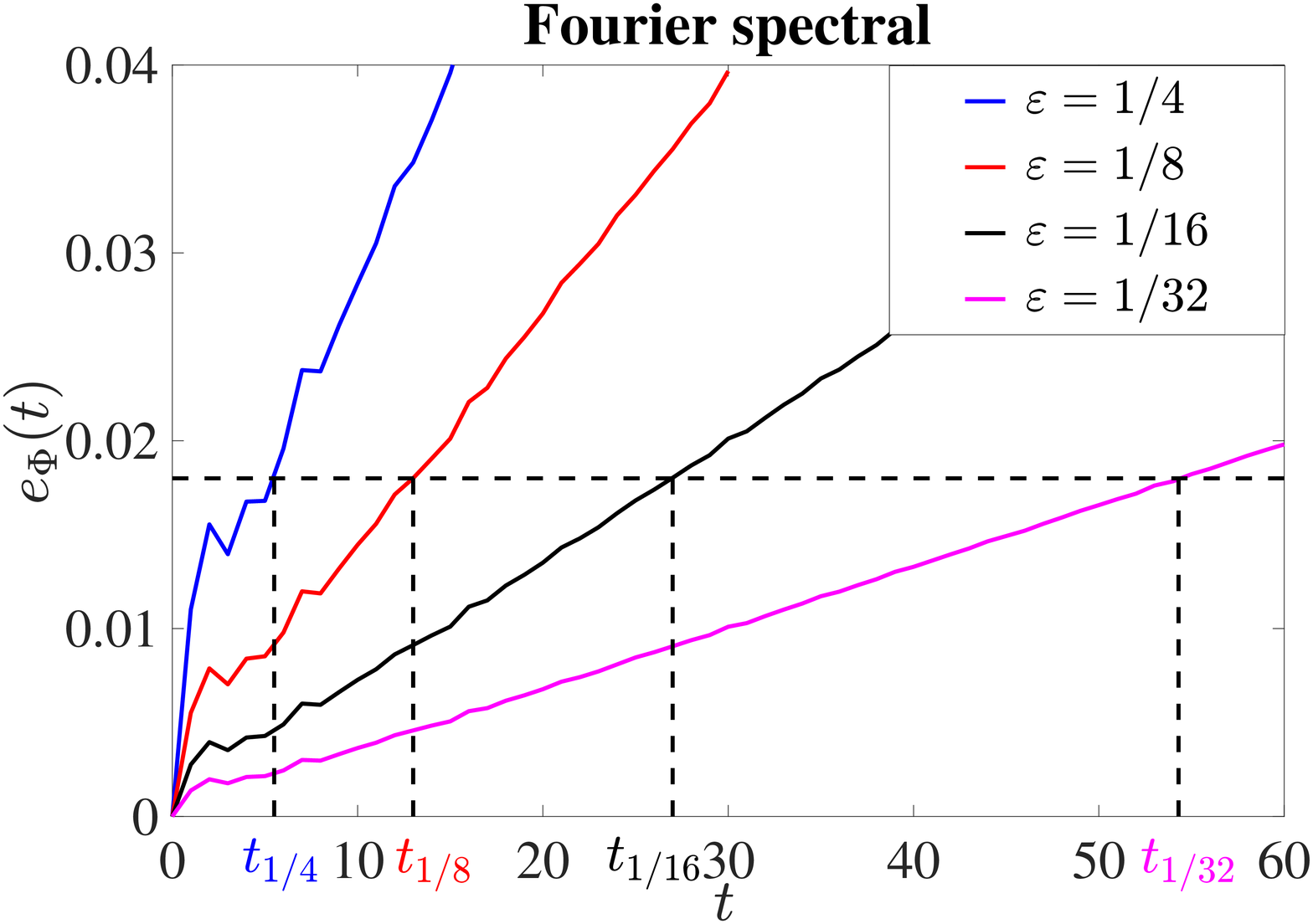}}
\end{minipage}
\caption{Spatial errors of finite difference and Fourier spectral spatial discretizations for the Dirac equation \eqref{eq:Dirac_21} in 1D with different $\varepsilon$.}
\label{fig:com_ep}
\end{figure}

Figure \ref{fig:comp} depicts long-time dynamics of the densities $\rho_j(2/\varepsilon, x) = |\phi_j(2/\varepsilon, x)|^2 (j =1, 2)$ of the Dirac equation in 1D with different $\varepsilon$ by the finite difference and spectral spatial discretizations. In the numerical simulations, we choose the mesh size as $h = \frac{\pi}{64}$ to compare the accuracy of these two spatial discretizations. For the fixed-time simulations, i.e., $\varepsilon = 1$, both spatial discretizations can get the numerical solution very accurately (cf. first row in Figure \ref{fig:comp}). For the long-time simulations, the spectral discretization can capture the long-time dynamics very accurately with the uniform mesh size, while the finite difference spatial discretization cannot get the correct numerical solution and the mesh size needs to be smaller when $\varepsilon$ becomes smaller (cf. middle and bottom rows in Figure \ref{fig:comp}).

Figure \ref{fig:com_ep} shows spatial errors of finite difference and Fourier spectral discretizations for the Dirac equation \eqref{eq:Dirac_21} in 1D with the initial data $\Phi_0(x) = (\cos(x), \cos(x))^T$, electromagnetic potentials $V(x) = 1/(1+\sin(x)^2)$ and $A_1(x) = (1+\sin(x))/(2+\cos(4x))$. The same mesh size $h = \pi/8$ and time step $\tau = 10^{-4}$ are taken for all numerical schemes. Results from the figure indicate the spatial errors of the Fourier spectral discretization are much smaller than the finite difference method in the long-time simulations. In addition, for the finite difference method, when $\eps$ is reduced by half, it needs twice the time to get $O(\varepsilon)$ spatial error, which confirms the spatial errors are at $O(1/\eps)$ in the long-time regime. Comparatively, for the Fourier spectral method, when $\eps$ is reduced by half, it needs twice the time to get the same spatial error, which validates the spatial errors are independent of $\eps$ in the long-time regime.

\subsection{Extension to the Dirac equation in the whole space}
In this subsection, we consider the following Dirac equation in the whole space for $d = 1, 2, 3$
\begin{equation}
\label{eq:Dirac_whole}
i\partial_t\Psi =  \left(- i\sum_{j = 1}^d
	\alpha_j\partial_j + \beta \right)\Psi+ \varepsilon\left(V(t, \mathbf{x})I_4 - \sum_{j = 1}^d A_j(t, \mathbf{x})\alpha_j\right)\Psi, \quad \mathbf{x}\in \mathbb{R}^d.
\end{equation}

By the same procedure, in 1D and 2D, the Dirac equation \eqref{eq:Dirac_whole} can be decoupled into two simplified PDEs with $\Phi := \Phi(t, \textbf{x}) = (\phi_1(t, \textbf{x}),\phi_2(t, \textbf{x}))^T \in \mathbb{C}^2$ satisfying 
\begin{equation}
i\partial_t \Phi = \left(-i\sum^d_{j=1} \sigma_j \partial_j +\sigma_3\right)\Phi	+ \left(V(t, \textbf{x})I_2 - \sum^d_{j=1}A_j(t,\textbf{x})\sigma_j\right)\Phi,
\label{eq:Dirac_whole21}
\end{equation}
for $x\in \mathbb{R}^d (d=1, 2)$, where $\Phi=(\psi_1,\psi_4)^T$ (or $\Phi=(\psi_2,\psi_3)^T$ ).

\begin{figure}[ht!]
\centerline{\includegraphics[width=14cm,height=7cm]{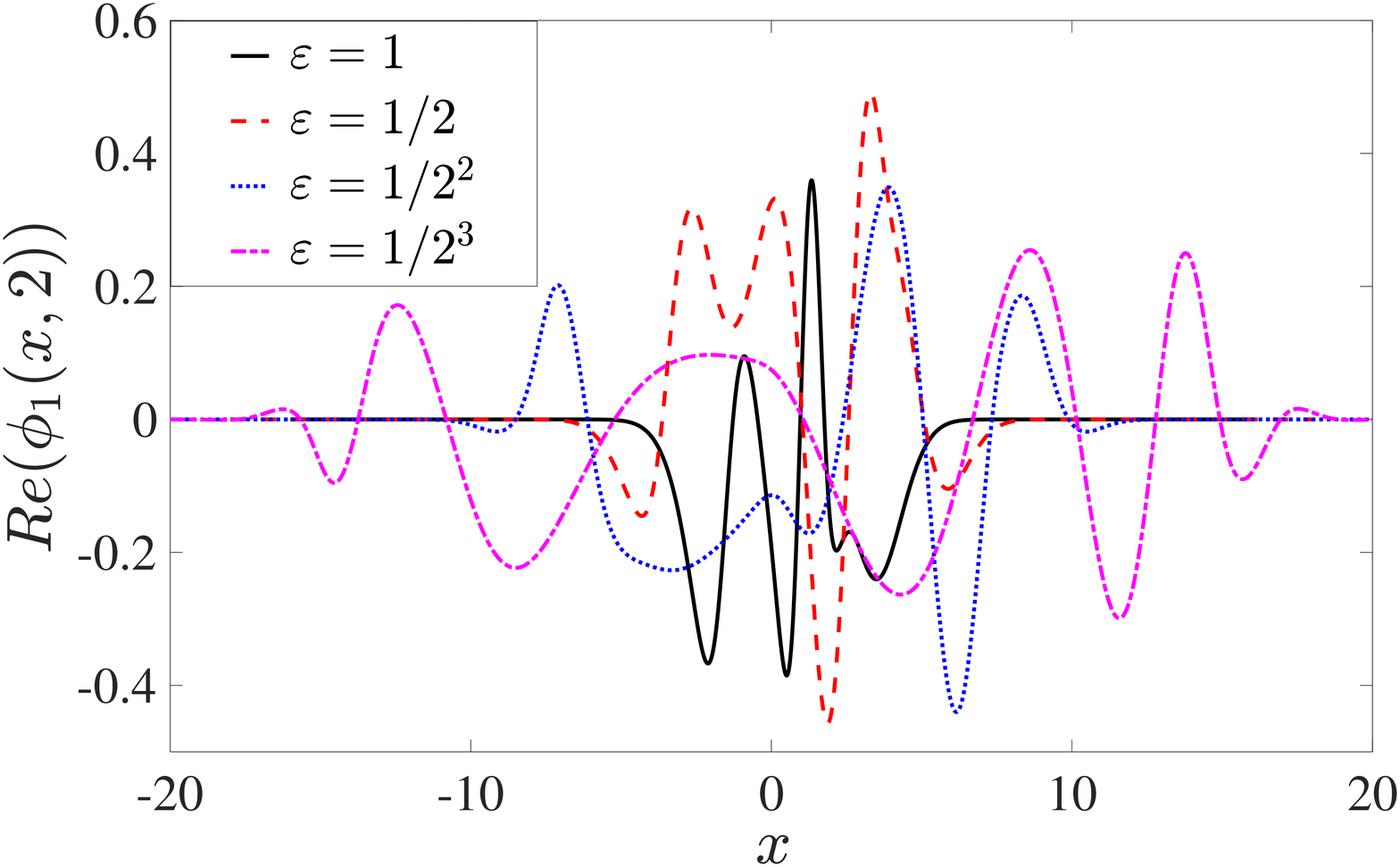}}
\caption{The solution $\phi_1(t=2/\varepsilon, x)$ of the Dirac equation \eqref{eq:Dirac_whole21} with $d= 1$ for different $\varepsilon$. Re($f$) denotes the real part of $f$.}
\label{fig:phi1_x}
\end{figure}

Similar to most works in the literatures for the analysis and computation of the Dirac equation \citep{BCJY,MY}, we truncate the whole space problem onto a large enough domain $\Omega$ such that the truncation error can be neglected, and assert periodic boundary conditions. In the long-time dynamics of the Dirac equation, the solution propagates the out-spreading waves, which is quite different from the solutions of the Dirac equation in the nonrelativistic limit regime or semiclassical regime. To illustrate this further, Figure \ref{fig:phi1_x} shows the solution of the Dirac equation \eqref{eq:Dirac_whole21} with $d = 1$, $V(t, x) = \frac{1-x}{1+x^2}$, $A_1(t, x) = \frac{(1+x)^2}{1+x^2}$ and $\Phi_0(x) = (\exp(-x^2/2), \exp(-(x-1)^2/2))^T$ for different $\varepsilon$. Due to the out-spreading waves, we need to choose the $\varepsilon$-dependent domain $\Omega_{\varepsilon} = (a- T_0/\varepsilon, b+ T_0/\varepsilon)$  in practical computations.

In the numerical simulations, we take $d=1$ and the electromagnetic potentials to be
\begin{equation}
V(t, x) = \frac{1-x}{1+x^2},\quad A_1(t, x)=\frac{(x+1)^2}{1+x^2},\quad x\in \Omega_\eps,
\end{equation}
we further take the initial data as \citep{MY,WHJY}
\begin{equation}
\begin{split}
&\phi_1(0, x) = \frac{1}{2}e^{-4x^2}e^{i S_0(x)}\left(1+\sqrt{1+S'_0(x)^2}\right),\\
&\phi_2(0, x) = \frac{1}{2}e^{-4x^2}e^{i S_0(x)}S'_0(x),\quad x \in \Omega_\eps,
\end{split}
\end{equation}
with 
\begin{equation}
S_0(x) = \frac{1}{40}\left(1+\cos(2\pi x)\right), \quad x \in \Omega_\eps.	
\end{equation}
The problem is solved numerically on the $\varepsilon$-dependent interval $\Omega_{\varepsilon} = (-7-T_0/\varepsilon, 7 + T_0/\varepsilon)$ with periodic boundary conditions. The `reference' solution $\Phi(t, x) = (\phi_1(t, x), \phi_2(t, x))^T$ is obtained numerically by the TSFP method with a very fine mesh size and a very small time step, i.e., $h_e = 1/ 512$ and $\tau_e = 10^{-4}$. The errors are displayed at $t = 1/\varepsilon$.

From Table \ref{tab:whole_CNFD_phi}, we can observe second order convergence in space and time for the CNFD method to solve the whole space problem (cf. first row in Table \ref{tab:whole_CNFD_phi}). The $\varepsilon$-resolution of the CNFD method is $h = O(\varepsilon^{1/2})$ and $\tau = O(\varepsilon^{1/2})$, which is verified through the upper triangle parts of the table above the bold diagonal line. Tables \ref{tab:whole_CNFD_rho} and \ref{tab:whole_CNFD_J} display the errors of the total probability $e_{\rho} (t = 1/\varepsilon)$ and the current density $e_{\mathbf{J}}(t = 1/\varepsilon)$ by the CNFD method, which coincide with our error analysis again. The numerical results of the other three FDTD methods are similar, and we omit them for brevity.

\begin{table}[h!]
\caption{Spatial and temporal errors of the CNFD method for the wave function $e_{\Phi}(t=1/\varepsilon)$ of the Dirac equation \eqref{eq:Dirac_whole21} in 1D.}
\label{tab:whole_CNFD_phi}
\renewcommand{\arraystretch}{1.3}
\centering
\vspace{0.3cm}
\begin{tabular}{cccccc}
\hline
$e_{\Phi}(t=1/\varepsilon)$ &$\begin{matrix} h_0 = \frac{1}{16} \\ \tau_0 = 0.05 \end{matrix}$ & $\begin{matrix} h_0/2 \\ \tau_0/2 \end{matrix}$ & $\begin{matrix} h_0/2^2 \\ \tau_0/2^2 \end{matrix}$ &  $\begin{matrix} h_0/2^3 \\ \tau_0/2^3 \end{matrix}$ &  $\begin{matrix} h_0/2^4 \\ \tau_0/2^4 \end{matrix}$   \\
\hline
$\varepsilon = 1$  & \bf{3.87E-2} & 1.01E-2 & 2.54E-3 & 6.36E-4 & 1.59E-4  \\
order & \bf{-} & 1.94 & 1.99 & 2.00 & 2.00  \\
\hline
$\varepsilon = 1/4$  & 1.04E-1 & \bf{2.95E-2} & 7.54E-3 & 1.89E-3 & 4.72E-4  \\
order & - & \bf{1.82} & 1.97 & 2.00 & 2.00  \\
\hline
$\varepsilon = 1/4^2$  & 2.84E-1 & 9.75E-2 & \bf{2.73E-2} & 6.95E-3 & 1.74E-3  \\
order & - & 1.54 & \bf{1.84} & 1.97 & 2.00 \\
\hline
$\varepsilon = 1/4^3$  & 5.74E-1 & 2.83E-1 & 9.66E-2 & \bf{2.70E-2} & 6.85E-3  \\
order & - & 1.02 & 1.55 & \bf{1.84} & 1.98 \\
\hline
$\varepsilon = 1/4^4$  & 8.03E-1 & 5.73E-1 & 2.83E-1 & 9.65E-2 & \bf{2.70E-2}  \\
order & - & 0.49 & 1.02 & 1.55 & \bf{1.84} \\
\hline
\end{tabular}
\end{table}

\begin{table}[h!]
\caption{Spatial and temporal errors of the CNFD method for the total probability $e_{\rho}(t=1/\varepsilon)$ of the Dirac equation \eqref{eq:Dirac_whole21} in 1D.}
\label{tab:whole_CNFD_rho}
\renewcommand{\arraystretch}{1.3}
\centering
\vspace{0.3cm}
\begin{tabular}{cccccc}
\hline	
$e_{\rho}(t=1/\varepsilon)$ &$\begin{matrix} h_0 = \frac{1}{16} \\ \tau_0 = 0.05 \end{matrix}$ & $\begin{matrix} h_0/2 \\ \tau_0/2 \end{matrix}$ & $\begin{matrix} h_0/2^2 \\ \tau_0/2^2 \end{matrix}$ &  $\begin{matrix} h_0/2^3 \\ \tau_0/2^3 \end{matrix}$ &  $\begin{matrix} h_0/2^4 \\ \tau_0/2^4 \end{matrix}$ \\
\hline
$\varepsilon = 1$  & \bf{3.39E-2} & 9.19E-3 & 2.35E-3 & 5.91E-4 & 1.48E-4  \\
order & \bf{-} & 1.88 & 1.97 & 1.99 & 2.00  \\
\hline
$\varepsilon = 1/4$  & 9.55E-2 & \bf{2.25E-2} & 5.37E-3 & 1.29E-3 & 3.22E-4  \\
order & - & \bf{2.09} & 2.07 & 2.06 & 2.00  \\
\hline
$\varepsilon = 1/4^2$  & 2.10E-1 & 6.13E-2 & \bf{1.40E-2} & 3.58E-3 & 9.03E-4  \\
order & - & 1.78 & \bf{2.13} & 1.97 & 1.99 \\
\hline
$\varepsilon = 1/4^3$  & 2.81E-1 & 1.26E-1 & 3.80E-2 & \bf{9.78E-3} & 2.44E-3  \\
order & - & 1.16 & 1.73 & \bf{1.96} & 2.00 \\
\hline
$\varepsilon = 1/4^4$  & 3.74E-1 & 2.32E-1 & 1.07E-1 & 3.32E-2 & \bf{8.69E-3} \\
order & - & 0.69 & 1.12 & 1.69 & \bf{1.93} \\
\hline
\end{tabular}
\end{table}

\begin{table}[h!]
\caption{Spatial and temporal errors of the CNFD method for the current density $e_{\mathbf{J}}(t=1/\varepsilon)$ of the Dirac equation \eqref{eq:Dirac_whole21} in 1D.}
\label{tab:whole_CNFD_J}
\renewcommand{\arraystretch}{1.3}
\centering
\vspace{0.3cm}
\begin{tabular}{cccccc}
\hline
$e_{\mathbf{J}}(t=1/\varepsilon)$ &$\begin{matrix} h_0 = \frac{1}{16} \\ \tau_0 = 0.05 \end{matrix}$ & $\begin{matrix} h_0/2 \\ \tau_0/2 \end{matrix}$ & $\begin{matrix} h_0/2^2 \\ \tau_0/2^2 \end{matrix}$ &  $\begin{matrix} h_0/2^3 \\ \tau_0/2^3 \end{matrix}$ &  $\begin{matrix} h_0/2^4 \\ \tau_0/2^4 \end{matrix}$ \\
\hline
$\varepsilon = 1$  & \bf{4.84E-2} & 1.28E-2 & 3.24E-3 & 8.13E-4 & 2.03E-4  \\
order & \bf{-} & 1.92 & 1.98 & 1.99 & 2.00  \\
\hline
$\varepsilon = 1/4$  & 1.15E-1 & \bf{2.84E-2} & 6.84E-3 & 1.69E-3 & 4.22E-4  \\
order & - & \bf{2.02} & 2.05 & 2.02 & 2.00  \\
\hline
$\varepsilon = 1/4^2$  & 2.63E-1 & 7.49E-2 & \bf{1.79E-2} & 4.55E-3 & 1.14E-3  \\
order & - & 1.81 & \bf{2.07} & 1.98 & 2.00 \\
\hline
$\varepsilon = 1/4^3$  & 3.79E-1 & 1.67E-1 & 5.00E-2 & \bf{1.28E-2} & 3.20E-3  \\
order & - & 1.18 & 1.74 & \bf{1.97} & 2.00 \\
\hline
$\varepsilon = 1/4^4$  & 5.47E-1 & 3.32E-1 & 1.47E-1 & 4.47E-2 & \bf{1.16E-2} \\
order & - & 0.72 & 1.18 & 1.72 & \bf{1.95} \\
\hline
\end{tabular}
\end{table}

\section{Conclusions}
The finite difference discretization in time combined with different spatial discretizations, including finite difference method and Fourier spectral method, were adapted to solve the Dirac equation with small electromagnetic potentials, while the strength of the potentials is characterized by $\varepsilon$ with $0 < \varepsilon \leq 1$ a dimensionless parameter. Performance of the FDTD and FDFS methods for the long-time dynamics of the Dirac equation up to the time at $O(1/\varepsilon)$ was compared based on rigorous error estimates and numerical results. The error bounds depend explicitly on the mesh size $h$, the time step $\tau$ as well as the small parameter $\varepsilon \in (0, 1]$, which indicate the spatial and temporal resolution capacities of the numerical methods for the long-time dynamics of the Dirac equation. For the FDTD methods, in order to get ``correct'' numerical approximations of the Dirac equation up to the time at $O(1/\varepsilon)$, the $\varepsilon$-scalability (or meshing strategy requirement) should be taken a: $h = O(\varepsilon^{1/2})$ and $\tau = O(\varepsilon^{1/2})$. For the FDFP methods, the spatial error is uniform for $\varepsilon \in (0, 1]$, which performs better than the FDTD methods in space, especially when $0 < \varepsilon \ll 1$. Extensive numerical results were reported to confirm our error bounds and demonstrate that they are sharp.

\section*{Acknowledgements}

The authors would like to thank Professor Weizhu Bao for his valuable suggestions and comments. This work was partially supported by the Ministry of Education of Singapore grant R-146-000-290-114. Part of the work was done when the authors were visiting the Institute for Mathematical Sciences at the National University of Singapore in 2020.


\begin{thebibliography}{10}
 
    \bibitem{AMP}
   {D.A. Abanin, S.V. Morozov, L.A. Ponomarenko, R.V. Gorbachev, A.S. Mayorov, M.I. Katsnelson, K. Watanabe, T. Taniguchi, K.S. Novoselov, L.S. Levitov, A.K. Geim}, 
	Giant nonlocality near the Dirac point in graphene, Science 332 (2011) 328--330.  
	
	\bibitem{AZ}
	{M.J. Ablowitz, Y. Zhu},
    Nonlinear waves in shallow honeycomb lattices, SIAM J. Appl. Math. 72 (2012) 240--260.  
      
   \bibitem{AH}
   {E. Ackad, M. Horbatsch},
   Numerical solution of the Dirac equation by a mapped Fourier grid method, J. Phys. A: Math. General 38 (2005) 3157--3171.  
	
    \bibitem{Anderson}
	{C.D., Anderson},
	The positive electron, Phys. Rev. 43 (1933) 491--498.      

   \bibitem{AS}
   {A. Arnold, H. Steinr\"uck},
   The `electromagnetic' Wigner equation for an electron with spin, ZAMP 40 (1989) 793--815.  

    \bibitem{BC}
	{W. Bao, Y. Cai},
	Mathematical theory and numerical methods for {B}ose-{E}instein condensation, Kinet. Relat. Models 6 (2013) 1-135.
	
     \bibitem{BCJT2016}
	{W. Bao, Y. Cai, X. Jia, Q. Tang},
	 A uniformly accurate multiscale time integrator pseudospectral method for the Dirac equation in the nonrelativistic limit regime, SIAM J. Numer. Anal. 54 (2016) 1785-1812.	

     \bibitem{BCJT}
	{W. Bao, Y. Cai, X. Jia, Q. Tang},
	Numerical methods and comparison for the Dirac equation in the nonrelativistic limit regime, J. Sci. Comput. 71 (2017) 1094-1134.
    	
     \bibitem{BCJY}
     {W. Bao, Y. Cai, X. Jia, J. Yin}, 
     Error estimates of numerical methods for the nonlinear Dirac equation in the nonrelativistic limit regime, Sci. China Math. 59 (2016) 1461-1494.
        
     \bibitem{BCY}
     {W. Bao, Y. Cai, J. Yin},
     Uniform error bounds of time-splitting methods for the nonlinear Dirac equation in the nonrelativistic limit regime, arXiv:1906.11101.
      
      \bibitem{BCY2}
      {W. Bao, Y. Cai, J. Yin},
      Super-resolution of time-splitting methods for the Dirac equation in the nonrelativistic regime, Math. Comp. 89 (2020) 2141-2173.
      
        \bibitem{BFY}
       {W. Bao, Y. Feng, W. Yi}, 
       Long time error analysis of finite difference time domain methods for the nonlinear Klein-Gordon equation with weak nonlinearity, Commun. Comput. Phys. 26 (2019) 1307-1334.
        
      \bibitem{BY}
     {W. Bao, J. Yin},
      A fourth-order compact time-splitting Fourier pseudospectral method for the Dirac equation, Res. Math. Sci. 6 (2019) article 11.
      
              
      \bibitem{BJP}
     {W. Bao, S. Jin, P.A. Markowich},
      On time-splitting spectral approximations for the Schr\"odinger equation in the semiclassical regime, J. Comput. Phys. 175 (2002) 487-524.
         
      \bibitem{BSG}
     {J.W., Braun, Q. Su, R. Grobe},
      Numerical approach to solve the time-dependent Dirac equation, Phys. Rev. A 59 (1999) 604--612.
         
       \bibitem{BHM}
	   {D. Brinkman, C. Heitzinger, P.A. Markowich}, 
	   A convergent 2D finite-difference scheme for the Dirac-Poisson system and the simulation of graphene, J. Comput. Phys. 257 (2014) 318--332.
	   
	   \bibitem{CMS}
	   {R. Carles, P.A. Markowich, C. Sparber}, 
	   Semiclassical asymptotics for weakly nonlinear {B}loch waves, J. Statist. Phys. 117 (2004) 343-375.
         
        \bibitem{Das1}
        {A. Das},
        General solutions of {M}axwell-{D}irac equations in 1 + 1-dimensional space-time and spatially confined solution, J. Math. Phys. 34 (1993) 3986--3999.
        
        \bibitem{Das2}
        {A. Das, D. Kay},
        A class of exact plane wave solutions of the {M}axwell-{D}irac equations, J. Math. Phys. 30 (1989) 2280--2284.
  
       \bibitem{Dirac1}
       {P.A.M. Dirac},
       The quantum theory of the electron, Proc. R. Soc. Lond. A 117 (1928) 610--624.
   
       \bibitem{Dirac2}
       {P.A.M. Dirac},
       The {P}rinciples of {Q}uantum {M}echanics, fourth ed., Oxford University Press, New York, 1958.
       
       \bibitem{Dujardin1}
       {G. Dujardin, E. Faou}, 
       Long time behavior of splitting methods applied to the linear Schr\"odinger equation, C. R. Math. Acad. Sci. Paris 344 (2007) 89--92.
       
       \bibitem{Dujardin2}
       {G. Dujardin, E. Faou},
       Normal form and long time analysis of splitting schemes for the linear Schr\"odinger equation with small potential, Numer. Math. 108  (2007) 223--262.
   
        \bibitem{ES}
       {M. Esteban, E. S\'er\'e}, 
        Existence and multiplicity of solutions for linear and nonlinear Dirac problems, Partial Differ. Equ. Appl. 12 (1997) 107--112.
                
        \bibitem{FW1}
        {C.L. Fefferman, M.I. Weinstein}, 
        Honeycomb lattice potentials and Dirac points, J. Am. Math. Socs 25 (2012) 1169--1220.
         
        \bibitem{FW2}
       {C.L. Fefferman, M.I. Weinstein}, 
        Wave packets in honeycomb structures and two-dimensional Dirac equations, Commun. Math. Phys. 326 (2014) 251--286.
        
        \bibitem{Feng}
        {Y. Feng},
        Long time error analysis of the fourth-order compact finite difference methods for the nonlinear Klein-Gordon equation with weak nonlinearity, Numer. Methods Partial Differential Equations 37 (2021) 897-914.
        
        \bibitem{FLB1}
        {F. Fillion-Gourdeau, E. Lorin, A.D. Bandrauk}, 
        Numerical solution of the time-dependent Dirac equation in coordinate space without fermion-doubling, Comput. Phys. Commun. 183 (2012) 1403--1415.
        
        \bibitem{FLB}
        {F. Fillion-Gourdeau, E. Lorin, A.D. Bandrauk}, 
        Resonantly enhanced pair production in a simple diatomic model, Phys. Rev. Lett. 110 (2013) 013002.
             
        \bibitem[Fillion-Gourdeau \emph{et al.}(2014)]{FLB2}
        {F. Fillion-Gourdeau, E. Lorin, A.D. Bandrauk}, 
        A split-step numerical method for the time-dependent Dirac equation in 3-D axisymmetric geometry, J. Comput. Phys. 272 (2014) 559--587.
        
        \bibitem{GGT}
        {F. Gesztesy, H. Grosse, B. Thaller}, 
        A rigorous approach to relativistic corrections of bound state energies
for spin-1/2 particles, Ann. Inst. Henri Poincar\'e Phys. Theor. 40 (1984) 159--174.

        \bibitem{Gosse}
        {L. Gosse},
         A well-balanced and asymptotic-preserving scheme for the one-dimensional linear Dirac equation, BIT Numer. Math. 55 (2015) 433--458.

        \bibitem{Gross}
        {L. Gross}, 
        The Cauchy problem for the coupled {M}axwell and {D}irac equations, Commun. Pure Appl. Math. 19 (1966) 1--15.
                     
        \bibitem{GSX}
        {B.-Y. Guo, J. Shen, C.-L. Xu} 
        Spectral and pseudospectral approximations using {H}ermite
functions: {A}pplication to the {D}irac equation, Adv. Comput. Math. 19 (2003) 35--55.
        
        \bibitem{JWY}
        {S. Jin, H. Wu, X. Yang}, 
        A numerical study of the Gaussian beam methods for one-dimensional Schr\"odinger-Poisson equations, J. Comput. Math. 28 (2010) 261--272.
            
        \bibitem{LMZ}
        {H. Liang, J. Meng, S-G. Zhou}, 
        Hidden pseudospin and spin symmetries and their origins in atomic nuclei, Phys. Rep. 570 (2015) 1--84.      
         
        \bibitem{MY}
        {Y. Ma, J. Yin}, 
         Error bounds of the finite difference time domain methods for the Dirac equation in the semiclassical regime, J. Sci. Comput. 81 (2019) 1801--1822.
         
        \bibitem{NGPNG}
        {A.H.C. Neto, F. Guinea, N.M.R. Peres, K.S. Novoselov, A.K. Geim},
         The electronic properties of graphene, Rev. Mod. Phys., 81 (2009) 109-162.
         
        \bibitem{NGMJ}
        {K.S. Novoselov, A.K. Geim, S.V. Morozov, D. Jiang, Y. Zhang, S.V. Dubonos, I.V. Grigorieva, A.A. Firsov}, 
        Electric field effect in atomically thin carbon films, Science 306 (2004) 666--669.
         
        \bibitem{Ring}
        {P. Ring}, 
         Relativistic mean field theory in finite nuclei, Prog. Part. Nucl. Phys. 37 (1996) 193--263.  
         
         
        \bibitem{Shebalin}
        {J.V. Shebalin},
        Numerical solution of the coupled Dirac and Maxwell equations, Phys. Lett. A 226  (1997) 1--16.
         
        \bibitem{ST}
	    {J. Shen, T. Tang}, 
	    Spectral and High-Order Methods with Applications, Science Press, Beijing, 2006.
         
        \bibitem{Smith}
       {G.D. Smith},
        Numerical Solution of Partial Differential Equations: Finite Difference Methods, Clarendon Press, Oxford, 1985.
        
        \bibitem{SM1}
       {C. Sparber, P.A. Markowich}, 
        Semiclassical asymptotics for the Maxwell-Dirac system, J. Math. Phys. 44 (2003) 4555--4572.
                  
         \bibitem{Spohn}
        {H. Spohn},
        Semiclassical limit of the Dirac equation and spin precession, Ann. Physics 282 (2000) 420--431.
    
         \bibitem{TE}
        {E. Tadmor},
        A review of numerical methods for nonlinear partial differential equations, Bull. Amer. Math. Soc. 49 (2012) 507--554.
        
        \bibitem{WHJY}
       {H. Wu, Z. Huang, S. Jin, D. Yin}, 
        Gaussian beam methods for the Dirac equation in the semi-classical regime, Commun. Math. Sci. 10 (2012) 1301--1315.
         
        \bibitem{XST}
        {J. Xu, S. Shao, H. Tang}, 
        Numerical methods for nonlinear Dirac equation, J. Comput. Phys. 245 (2013) 131--149.
        

\end{thebibliography}
\end{document}